\documentclass[onefignum,onetabnum]{siamart250211}

\usepackage{lipsum}
\usepackage{amsfonts}
\usepackage{graphicx}
\usepackage{epstopdf}
\usepackage[noend]{algorithmic}
\usepackage{subcaption}
\usepackage{tabularray}
\usepackage{bm}
\usepackage{multirow}
\ifpdf
  \DeclareGraphicsExtensions{.eps,.pdf,.png,.jpg}
\else
  \DeclareGraphicsExtensions{.eps}
\fi

\newsiamremark{remark}{Remark}
\newsiamremark{hypothesis}{Hypothesis}
\crefname{hypothesis}{Hypothesis}{Hypotheses}
\newsiamthm{claim}{Claim}

\headers{A Cubed Sphere Fast Multipole Method}{A. Chen and R. Krasny}

\title{A Cubed Sphere Fast Multipole Method\thanks{Submitted to the editors \today.
\funding{The first author was supported by 
NSF Graduate Research Fellowship grant DGE-2241144
and
the second author was supported by
NSF grant DMS-2110767.}}}

\author{Anthony Chen\thanks{Department of Mathematics, University of Michigan, Ann Arbor, MI 
  (\email{cygnari@umich.edu}).}
\and Robert Krasny\thanks{Department of Mathematics, University of Michigan, Ann Arbor, MI
  (\email{krasny@umich.edu}).}}

\usepackage{amsopn}

\makeatletter
\newcommand*{\addFileDependency}[1]{
  \typeout{(#1)}
  \@addtofilelist{#1}
  \IfFileExists{#1}{}{\typeout{No file #1.}}
}
\makeatother

\ifpdf
\hypersetup{
  pdftitle={A Cubed Sphere Fast Multipole Method},
  pdfauthor={A. Chen and R. Krasny}
}
\fi

\begin{document}

\maketitle

\begin{abstract}
This work describes a new version of the
Fast Multipole Method for 
summing pairwise particle interactions 
that arise from discretizing
integral transforms and
convolutions on the sphere.
The kernel approximations use
barycentric Lagrange interpolation on
a quadtree composed of cubed sphere grid cells.
The scheme is kernel-independent
and
requires kernel evaluations 
only at points on the sphere.
Results are presented for the
Poisson and biharmonic equations on the sphere,
barotropic vorticity equation on a rotating sphere,
and
self-attraction and loading potential in
tidal calculations.
A tree code version is also described for
comparison,
and
both schemes are tested in
serial and parallel calculations. 
\end{abstract}



\section{Introduction}
\label{sec:intro}

\begin{figure}[htb]
\centering
\includegraphics[width=\linewidth]{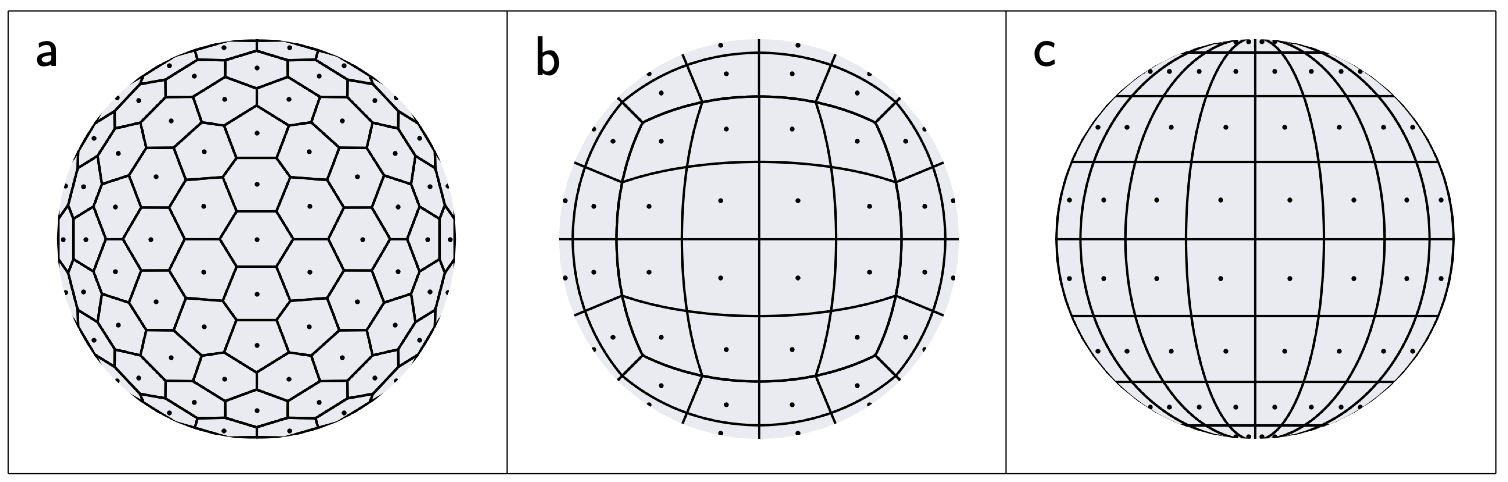}
\caption{Partitions of the sphere
showing grid cells and cell centers,
(a) icosahedral grid,
hexagonal cells
(aside from 12~pentagons),
cell centers are vertices of the
icosahedral triangulation,
(b) cubed sphere grid,
cell centers are mapped to the sphere
from uniform grid points on the faces 
of the inscribed cube, 
(c) latitude-longitude grid,
cell center coordinates are given by averaging 
the latitude and longitude of cell edges.}
\label{fig:grids}
\end{figure}

This work is concerned with computing
integral transforms
on the sphere $S \subset \mathbb{R}^3$,
\begin{equation}
\label{eq:integral_transform}
\phi({\bf x}) = 
\int_{S}K(\mathbf{x},\mathbf{y})
f(\mathbf{y})dS(\mathbf{y}),
\quad {\bf x} \in S,
\end{equation} 
where $f$ is an input data field,
the kernel $K$ encodes some process,
and
$\phi$ is the output field.
When the kernel is a Green's function,
it will be written $G$ instead of $K$.
Several methods are available
to compute the integral 
numerically~\cite{atkinson1982numerical,bellet2022quadrature,hesse2010numerical,lebedev1976quadratures,womersley2018efficient},
but for simplicity we consider a scheme in which
the sphere is partitioned into grid cells
as shown for three examples in
\Cref{fig:grids},
and~\cref{eq:integral_transform} is
approximated by the midpoint rule, 
\begin{equation}\label{eq:nbodysum}  \phi(\mathbf{x}_i) \approx
\sum_{j=1}^N K(\mathbf{x}_i,\mathbf{x}_j)w_j,
\quad
w_j = f(\mathbf{x}_j)A_j,\quad i=1:N,
\end{equation}
where the quadrature points ${\bf x}_i \in S$  
are the cell centers
and
$A_j$ are the geodesic cell areas.
We shall also view~\cref{eq:nbodysum} 
as a pairwise particle $N$-body sum,
where
${\bf x}_i$ is a target particle
and
${\bf x}_j$ is a source particle.
Direct summation of \cref{eq:nbodysum}
using loops over indices $i$ and $j$
requires $O(N^2)$ operations,
which is prohibitively expensive for large $N$.
Here we describe a new version of the
Fast Multipole Method~\cite{fmmgreengard}
for computing~\cref{eq:nbodysum}
with general kernels 
and
point sets on the sphere.
The scheme can accommodate
cases in which the target and source
particles are different sets,
but they are assumed to lie on the sphere.
If the kernel is singular for
${\bf x} = {\bf y}$,
then the $i=j$ term in the sum is omitted;
more accurate quadrature rules
can be handled as long as they have
the form~\eqref{eq:nbodysum},
perhaps with a local correction.

\subsection{Motivation}
\label{sec:motivation}

A number of applications in
astrophysics
and
geophysics
involve processing a data field 
$f$ defined on a sphere (see
for example~\cite{
brun2009large,parker1994geophysical,
suo2023spherical,
toth2011obtaining}).
A common approach first expands
the data field in spherical harmonics, 
\begin{equation}
\label{eq:SHT}
f(\theta,\lambda) =
\sum_{n=0}^{\infty}\sum_{m=-n}^n
\widehat{f}_n^mY_n^m(\theta,\lambda),
\quad \widehat{f}_n^m =
\int_0^{2\pi}\!\!\int_0^{\pi}
f(\theta,\lambda)\overline{Y_n^m}(\theta,\lambda)\sin\theta d\theta d\lambda,
\end{equation}
where 
$\widehat{f}_n^m$ is the 
$(n,m)$-spherical harmonic coefficient of $f$,
and
points on the sphere are expressed as 
${\bf x} = (\theta,\lambda)$ with
co-latitude~$\theta$ and longitude~$\lambda$.
The data coefficients are then multiplied by
coefficients $\alpha(n)$ 
that encode the process under consideration,
and
the modified coefficients are
transformed back to the spatial domain yielding the 
output field,
\begin{equation}
\label{eq:SH}
\phi(\theta,\lambda) =
\sum_{n=0}^{\infty}\sum_{m=-n}^n
\alpha(n) \widehat{f}_n^m
Y_n^m(\theta,\lambda).
\end{equation}
When this approach is used to solve a 
partial differential equation on the sphere, 
$\alpha(n)$ is the reciprocal of the 
operator symbol~\cite{boyd2001chebyshev}.

In practice,
the integrals in~\cref{eq:SHT} are typically 
discretized on a latitude-longitude grid
and
the infinite sums are truncated,
yielding finite spherical harmonics 
transforms (SHT). 
Fast methods 
for computing the SHT have been developed
(see for example~\cite{healy2003ffts,
mohlenkamp1999fast,
reinecke2013libsharp,
schaeffer2013efficient,
suda2002fast,
swarztrauber2000generalized,
tygert2010fast}), 
yet there are continuing challenges
including
the need to interpolate 
when the data field
is not given on a latitude-longitude grid, 
the high computational cost 
if many harmonics are needed
to resolve localized features,
the introduction of spurious oscillations 
when the infinite sums are truncated,
and the need for more efficient 
parallel implementations 
to handle large-scale problems.
Our objective here
is to explore an alternative 
to the spherical harmonics transform~\cref{eq:SH} 
for several problems in
geophysical fluid dynamics.

\subsection{Spherical convolution}

An alternative approach replaces the
spherical harmonics transform~\cref{eq:SH}
by a spherical convolution,
\begin{equation}
\label{eq:spherical_convolution}
\phi({\bf x}) = 
(f \ast g)({\bf x}) \equiv
\int_S g({\bf x}\cdot{\bf y})f({\bf y})dS({\bf y}),
\quad {\bf x}, {\bf y} \in S,
\end{equation}
where the connection to~\cref{eq:SH} 
is made through the spherical convolution
theorem~\cite{driscoll1994computing},
\begin{equation}
\label{eq:convolution_theorem}
(\widehat{f\ast g})_n^m =
\sqrt{\frac{4\pi}{2n+1}}\widehat{f}_n^m\widehat{g}_n^0.
\end{equation}
It follows that the convolution kernel is
composed of zonal spherical harmonics,
\begin{equation}
g({\bf x}\cdot{\bf y}) =
\sum_{n=0}^{\infty}
\widehat{g}_n^0 Y_n^0(\varphi,0),
\quad
\widehat{g}_n^0 =
\sqrt{\frac{2n+1}{4\pi}}\alpha(n),
\quad {\bf x}, {\bf y} \in S,
\end{equation}
where ${\bf x}\cdot{\bf y} = \cos\varphi$
and
$\varphi$ is the great circle angle between 
$\mathbf{x}$ and $\mathbf{y}$.
After noting that
$Y_n^0(\varphi,0) = 
\sqrt{(2n+1)/4\pi}P_n(\cos\varphi)$,
where $P_n$ is the
Legendre polynomial of degree $n$,
the spherical convolution~\cref{eq:spherical_convolution}
becomes an integral transform~\cref{eq:integral_transform},
where the kernel is given by the Legendre series
\begin{equation}
\label{eq:kernel_Legendre}
K(\mathbf{x},\mathbf{y}) =
\sum_{n=0}^{\infty}
\frac{2n+1}{4\pi}\alpha(n) 
P_n({\bf x}\cdot{\bf y}),
\quad \mathbf{x}, \mathbf{y} \in S.
\end{equation}
In some cases treated below
the process coefficients~$\alpha(n)$
and kernel $K$
are known in closed form,
but we shall also consider 
an example from tidal calculations
where the $\alpha(n)$
are given by empirical data
and
$K$ is only known approximately.
However, even assuming $K$ can be
computed accurately and efficiently,
the cost of computing the
$N$-body sum~\cref{eq:nbodysum}
must be addressed.

\subsection{Fast summation methods}

For the case of the Newtonian potential
and
general point sets in~$\mathbb{R}^3$,
several fast summation methods 
were developed that reduce 
the operation count for~\cref{eq:nbodysum}
to $O(N\log{N})$ or $O(N)$. 
Early methods include a dual tree traversal algorithm~\cite{appel1985efficient} 
and a tree code~\cite{barneshut} 
using far-field monopole approximations
of the kernel,
while the
Fast Multipole Method (FMM)~\cite{cheng1999fast,fmmgreengard} 
used higher order near-field and far-field
multipole approximations 
for enhanced accuracy and efficiency.
More recent developments include
kernel-independent methods using
equivalent densities~\cite{ying2004kernel},
polynomial interpolation~\cite{fong2009black,
wang2020kernel},
and
interpolative decomposition~\cite{martinsson2011randomized,
xing2020interpolative}.

For point sets on the sphere,
the simplest approach is to use a
fast summation method
for general point sets in $\mathbb{R}^3$.
For example,
a 3D Cartesian Taylor series tree code 
was used to compute particle interactions
on the sphere
given by the spherical Biot-Savart kernel~\cite{draghicescu1994efficient,sakajo2009extension}.
However, this required 
kernel evaluations at points off the sphere,
and there may be an advantage in 
designing methods that require
kernel evaluations only on the sphere;
such methods include an FMM-like method based on
interpolative decomposition for
the spherical Neumann Green's 
function~\cite{kaye2020fast},
and an FMM based on analytic series
approximations for the
spherical Laplace Green's 
function~\cite{suo2023spherical}.

\subsection{Present work}

Previous work developed a
tree code~\cite{bltc,wang2020kernel}
and
FMM~\cite{bldtt}
for kernels and point sets in~$\mathbb{R}^3$
using barycentric Lagrange interpolation
on an octree composed of rectangular 
boxes~\cite{berrut2004barycentric}.
The present work adapts this approach to 
kernels and point sets on the sphere
using 
barycentric Lagrange interpolation on 
a quadtree composed of cubed sphere grid cells.
The resulting schemes are called the
Cubed Sphere FMM (CSFMM) and
Cubed Sphere Tree Code (CSTC);
they are kernel-independent
and
require kernel evaluations
only at points on the sphere.
It should be noted that
the cubed sphere is used here
in two different ways, 
(1) it is one of the three 
spherical partitions shown in Figure 1
used to discretize the convolution integral,
(2) the cubed sphere grid cells provide
the particle clusters used in the 
CSFMM and CSTC to
accelerate the calculation of the $N$-body sum.

The rest of the article is organized as follows. 
\Cref{sec:sphereconv}
describes four problems whose solution
can be expressed as an integral
transform~\cref{eq:integral_transform};
these are the
Poisson and biharmonic equations on the sphere,
barotropic vorticity equation (BVE)
on a rotating sphere,
and
self-attraction and loading (SAL) potential
in tidal calculations.
\Cref{sec:cs} reviews the cubed sphere 
and
explains the tree building process.
\Cref{sec:bli} reviews barycentric Lagrange interpolation and its
extension to cubed sphere grid cells.
\Cref{sec:interactions} explains the particle
interactions used in the proposed methods.
\Cref{sec:upward} describes the upward pass
used in both methods.
\Cref{sec:cstc} and \Cref{sec:csfmm} 
present the CSTC
and 
CSFMM algorithms.
\Cref{sec:details} gives
implementation details.
\Cref{sec:error}
assesses numerical errors
and
\Cref{sec:runtime} gives the
runtime in serial and parallel calculations.
\Cref{sec:BVEres} and
\Cref{sec:SAL} present
BVE and SAL calculations.
A summary is given in~\Cref{sec:summary}.

\section{Problems}
\label{sec:sphereconv}

This section describes four problems 
whose solution can be expressed as an
integral 
transform~\eqref{eq:integral_transform}
for a certain kernel. 
The problems are posed on the unit sphere 
and spheres with non-unit radius can be handled 
by suitable scaling. 

\subsection{Poisson equation}
\label{sec:poisson}

The first problem is the Poisson equation
on the sphere,
\begin{equation}
\label{eq:Poisson_equation}
-\Delta\phi({\bf x}) = f({\bf x}),
\quad {\bf x} \in S,
\end{equation}
with the restriction that the data field has mean zero, 
$\int_S fdS = 0$.
An application appearing below
arises when $f$ is the vorticity of a fluid flow
on the sphere and
$\phi$ is the stream function.
The kernel is the spherical Laplace Green's function~\cite{bogomolov1977dynamics, kimura1987vortex},
\begin{equation}
\label{eq:G_L}
G_{\rm L}(\mathbf{x},\mathbf{y}) =
-\frac{1}{4\pi}\ln(1-\mathbf{x}\cdot\mathbf{y}),
\quad {\bf x}, {\bf y} \in S.
\end{equation}
Spherical harmonics are eigenfunctions 
of the Laplacian, 
$\Delta Y_n^m=-n(n+1)Y_n^m$,
and
this provides a reference
to test the accuracy of the proposed methods.
In practice we use the 
real spherical harmonics
$Y_{n,m} = 
(Y_n^m + (-1)^mY_n^{-m})/\sqrt{2}$ 
with order $m > 0$,
where the solution of the Poisson equation 
is written as a convolution,
\begin{equation}
\label{eq:G_L_1}
Y_{n.m}({\bf x}) =
-n(n+1)
\int_S G_L({\bf x},{\bf y})Y_{n,m}({\bf y})
dS({\bf y}),
\quad n \ge 1, 
\quad {\bf x} \in S.
\end{equation}
When the integral is discretized,
the computed $Y_{n,m}$ 
can be compared with the exact value.
For degree $n=0$ 
we have the constant 
$Y_{0,0} =1/(2\sqrt{\pi})$
and the following relation is used
as the reference,
\begin{equation}
\label{eq:G_L_2}
\int_S G_L({\bf x},{\bf y})Y_{0,0}({\bf y})
dS({\bf y}) = \frac{1 - \ln 2}{2\sqrt{\pi}}.
\end{equation}

\subsection{Biharmonic equation}

The biharmonic equation on the sphere,
\begin{equation}
\Delta^2\phi({\bf x}) = f({\bf x}),
\quad {\bf x} \in S,
\end{equation}
arises in geodesy and remote 
sensing~\cite{hardy1990theory,
sandwell1987biharmonic}. 
The kernel is the
spherical biharmonic
Green's function,
\begin{equation}
\label{eq:G_B}
G_{\rm B}(\mathbf{x},\mathbf{y}) =
\frac{1}{4\pi}\mathrm{dilog}
\left(\frac{1+\mathbf{x}\cdot\mathbf{y}}{2}\right),
\quad {\bf x}, {\bf y} \in S,
\end{equation}
where the dilogarithm is defined by 
\begin{equation}
\mathrm{dilog}(x) =
-\int_0^x\frac{\log(1-t)}{t}dt, \quad x \le 1.
\end{equation}
In this work the dilogarithm is computed
using a rational  
approximation~\cite{voigt2022comparison}. 
Various expressions for $G_{\rm B}$ 
appear in the
literature~\cite{freeden2008spherical,
parker1994geophysical},
and~\cref{eq:G_B} 
is consistent with the results in~\cite{wessel2008interpolation}. 
Spherical harmonics are eigenfunctions
of the biharmonic operator,
$\Delta^2Y_n^m=n^2(n+1)^2Y_n^m$,
again providing a reference solution.

\subsection{Barotropic vorticity equation}
\label{sec:BVE}

The next example is the 
barotropic vorticity equation (BVE) 
for incompressible inviscid fluid flow on a 
rotating sphere~\cite{vallis2017atmospheric}.
The Coriolis parameter is 
$f({\bf x}) = 2\Omega z$ with $z=\cos\theta$,
where the sphere rotates with
angular frequency $\Omega$ about the $z$-axis
in a Cartesian coordinate system.
The BVE expresses the 
conservation of absolute vorticity 
following fluid particles,
\begin{equation}
\label{eq:BVE}
\frac{D(\zeta+f)}{Dt} \equiv
\frac{\partial(\zeta+f)}{\partial t} +(\mathbf{u}\cdot\nabla)(\zeta+f) = 0,
\end{equation}
where 
$D/Dt$ is the material derivative,
$\zeta({\bf x},t)$ is the relative vorticity,  
and
$\mathbf{u}({\bf x},t)$ is the fluid velocity. 
An initial vorticity is given 
and the goal
is to study its evolution in time. 
We use a vortex method that tracks a set of
Lagrangian particles ${\bf x}_i(t)$
on the sphere with vorticity 
$\zeta_i(t)$~\cite{chorin1973numerical,
bosler2014lagrangian,
newton2007n}. 
This yields a system of ODEs for the
particle motion and 
the change in their vorticity,
\begin{subequations}
\label{eq:BVE_ODEs}
\begin{align}
\label{eq:K_BS}
\frac{d{\bf x}_i}{dt} &=
\sum_{j=1}^N
K_{\mathrm{BS}}(\mathbf{x}_i,\mathbf{x}_j)
\zeta_jA_j,
\quad
K_{\rm BS}({\bf x},{\bf y}) =
-\frac{1}{4\pi}
\frac{\mathbf{x} \times \mathbf{y}}
{1-\mathbf{x}\cdot\mathbf{y}},
\quad {\bf x},{\bf y} \in S, \\
\quad \frac{d\zeta_i}{dt} &=
-2\Omega \frac{dz_i}{dt},
\quad i=1:N,
\end{align}
\end{subequations}
where the spherical Biot-Savart kernel,
$K_{\mathrm{BS}}(\mathbf{x},\mathbf{y}) =
\nabla_{\mathbf{x}}
G_{\rm L}(\mathbf{x},\mathbf{y}) 
\times \mathbf{x}$,
is obtained from the
spherical Laplace Green's 
function~\cref{eq:G_L}.
The particles initially lie at the centers
of spherical grid cells with area $A_j$  
as in~\cref{fig:grids}.
The velocity calculation~\cref{eq:K_BS}
has the form of an $N$-body 
sum~\cref{eq:nbodysum}.
Remeshing is carried out at each 
time step~\cite{bosler2014lagrangian,koumoutsakos1997inviscid}.

\subsection{Self-attraction and loading}
\label{sec:sal} 

The accurate prediction of ocean tides 
is an
important factor in commercial shipping and
naval operations.
For such calculations the
shallow water equations (SWE)
can be supplemented by a forcing term that 
accounts for the
gravitational self-attraction of the water mass 
and elastic deformation of the Earth surface
due to the varying sea surface height~\cite{kuhlmann2015self}.
The forcing term is the
gradient of the self-attraction and loading (SAL)
potential
$\eta_{\rm SAL}({\bf x},t)$,
which depends on the
sea surface height anomaly 
$\eta({\bf x},t)$.
Here we are concerned with computing the
SAL potential,
while its application in tidal SWE calculations is 
reserved for future work.

Recent calculations of $\eta_{\rm SAL}$
use the spherical harmonics approach
described in \Cref{sec:motivation}~\cite{
barton2022global, 
brus2023scalable, 
shihora2022self}. 
The input data field is the
sea surface height 
anomaly~$\eta$,
and
the spectral multiplier coefficients are 
\begin{equation}\label{eq:salsh}
\alpha(n) =
\frac{\rho_w}{\rho_e}
\frac{3(1+k_n'-h_n')}{2n+1},
\end{equation}
where
$\rho_w/\rho_e$ is the ratio of
the mean mass density of 
seawater to that of the Earth, 
and 
$k_n', h_n'$ are the modified 
Load Love Numbers (LLNs)
that encode the material properties of the
Earth~\cite{farrell1972deformation,hendershott1972effects}. 
The LLNs are not known analytically
and 
in practice they are computed numerically
by solving an elasticity problem
that accounts for the 
Earth's density profile and seismic properties;
a commonly used dataset 
based on the elastic Earth model PREM 
provides LLNs up to 
$n_{\rm max} = 46343$~\cite{wang2012load}.
However,
this approach faces challenges
in resolving local features
due to high calculation cost
and
the limitations of current software.
In one example~\cite{barton2022global},
researchers used the SHTns package~\cite{schaeffer2013efficient}
which implements the SHT 
on a structured Gaussian mesh,
yet can only run on a single node,
while another recent 
work~\cite{brus2023scalable} 
computed the SHT directly on an 
unstructured 
mesh since it has an
efficient multinode 
parallelization.
Even so, cost was still a significant issue
in~\cite{barton2022global,brus2023scalable}
to the extent that
the SAL term was computed only at 
prescribed time intervals 
rather than every time step.

With the goal of improving the accuracy 
and efficiency of SAL calculations, 
we compute~$\eta_{\rm SAL}$ as a spherical convolution~\cref{eq:spherical_convolution},
where following~\cref{eq:kernel_Legendre},
the SAL kernel is
expressed as a Legendre series with LLN coefficients~\cite{farrell1972deformation,
hendershott1972effects,
kuhlmann2015self},
\begin{equation}   
\label{eq:G_SAL}
G_{\rm SAL}(\mathbf{x},\mathbf{y}) =
\frac{3\rho_w}{4\pi\rho_e}
\sum_{n=0}^{\infty}(1+k_n'-h_n')
P_n(\mathbf{x}\cdot\mathbf{y}),
\quad {\bf x}, {\bf y} \in S.
\end{equation}
However, instead of
truncating the series~\cref{eq:G_SAL},
which would introduce spurious oscillations,
we use a closed form approximation for 
$G_{\rm SAL}$ that seeks to incorporate 
the asymptotic behavior of the LLNs.
We found empirically that the
LLNs in the PREM dataset~\cite{wang2012load} 
can be approximated by 
\begin{equation}
k_n'\approx a_1/n, \quad 
h_n'\approx b_0+b_1/n,
\end{equation}
with coefficients 
$a_1=-2.7$, $b_0=-6.21196$, $b_1=6.1$.
Then on writing
\begin{equation}
1+k_n'-h_n' = 
(1-b_0) + (k_n'-(h_n'-b_0)) \approx
(1-b_0) + (a_1-b_1)/n,
\end{equation}
in~\cref{eq:G_SAL}, and using
identities from Section 8.921 and 8.926 of~\cite{gradshteyn2014table},
\begin{equation}
\label{eq:Legendre}
\sum_{n=0}^{\infty}P_n(x) =
\frac{1}{\sqrt{2(1-x)}},
\quad
\sum_{n=1}^{\infty}\frac{1}{n}P_n(x)=-\ln\left(\sqrt{\frac{1-x}{2}}+\frac{1-x}{2}\right),
\end{equation}
the SAL kernel~\eqref{eq:G_SAL} 
has the closed form approximation 
\begin{equation}
\label{eq:G_SAL_approx}
G_{\mathrm{SAL}}
(\mathbf{x},\mathbf{y}) \approx
\frac{3\rho_w}{4\pi\rho_e}
\left(\frac{1-b_0}
{\gamma({\bf x},{\bf y})} -
(a_1-b_1)
\ln\left(\frac{\gamma({\bf x},{\bf y})}{2}
+ \frac{\gamma^2({\bf x},{\bf y})}{4}\right)
\right),
\end{equation}
where
$\gamma(\mathbf{x},\mathbf{y}) =
\sqrt{2(1-\mathbf{x}\cdot\mathbf{y})}$.
The SAL potential $\eta_{\rm SAL}$
can then be computed as a spherical convolution~\cref{eq:spherical_convolution}
of the $G_{\rm SAL}$ 
approximation~\cref{eq:G_SAL_approx}
with the sea surface height anomaly~$\eta$,
and
upon discretization this again leads to an
$N$-body sum~\cref{eq:nbodysum}.

\section{Cubed sphere and tree building}
\label{sec:cs}
\begin{figure}[htb]
\includegraphics[width=\textwidth]{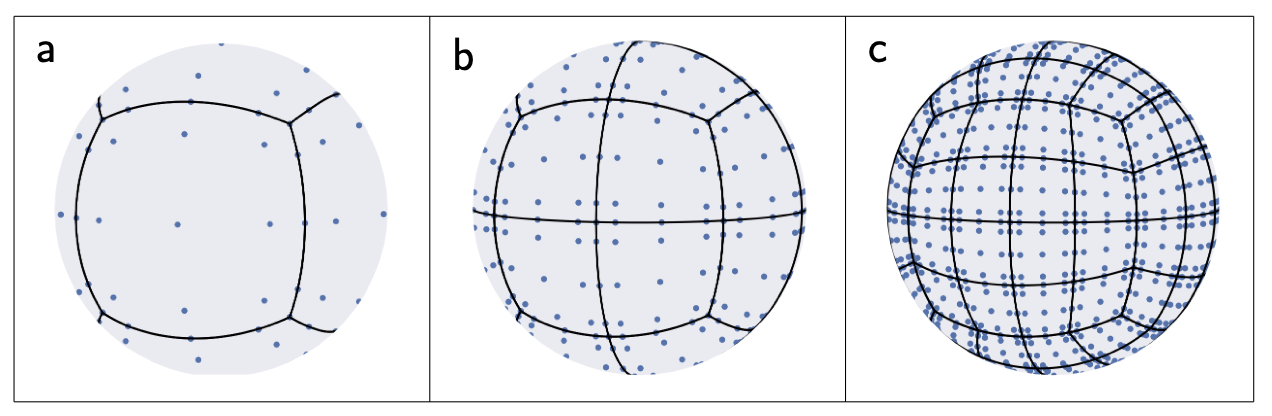}
\vskip -5pt
\caption{
Quadtree composed of
equiangular gnomonic cubed sphere grid cells,
(a) level~0,
(b) level~1,
(c) level~2,
CSFMM and CSTC use
barycentric Lagrange interpolation in cells,
example of $5 \times 5$ 
Chebyshev points~(\textcolor{blue}{$\bullet$})
in each cell.}
\label{fig:eacs}
\end{figure}

The cubed sphere is obtained by
projecting grid points 
from the faces of a cube 
to the circumscribed sphere.
In this work the cubed sphere is used in two ways,
(1) as one of the 
three spherical partitions
shown in~\Cref{fig:grids}
used to discretize the integral
transform~\cref{eq:integral_transform},
(2) the CSFMM and CSTC
use a quadtree
composed of cubed sphere grid cells 
with several levels of refinement,
where each cell
defines a cluster of particles,
and
the kernel approximations are done by
barycentric Lagrange interpolation
in the cells.
Among several versions of the cubed sphere, 
we use the equiangular gnomonic 
version~\cite{ronchi1996cubed},
as shown in \cref{fig:eacs}
with three levels of refinement.
This version 
has several favorable properties 
for interpolation and quadrature;
the cell edges lie on great circles,
the cells have nearly uniform area
(the area ratio between the smallest and largest cells is approximately $0.7$),
and
each cell can be parameterized as a
square in terms of 
local angle coordinates
$(\xi,\eta)$~\cite{ronchi1996cubed}.

To build the tree, 
start by assigning each particle ${\bf x}_i$
to one of the six faces of the 
cubed sphere
which define the root clusters of the tree.
If a cluster has more than a prescribed
number of particles~$N_0$, 
it is partitioned into four subclusters, 
and
this continues recursively
until reaching
the leaf clusters which
contain fewer than $N_0$ particles
(the choice of $N_0$ will be specified below).
This process yields a quadtree of 
particle clusters on the sphere,
where each cluster corresponds to a
cubed sphere grid cell.
When the particles are nonuniformly distributed,
it is advantageous to shrink the 
clusters to tightly bound their 
particles~\cite{li2009cartesian}.
The next section reviews
barycentric Lagrange 
interpolation~\cite{berrut2004barycentric}
and explains its application to 
cubed sphere grid cells.

\section{Barycentric Lagrange interpolation}
\label{sec:bli}

Consider a function $f$ defined on $[-1,1]$
and
let $s_k = \cos(k\pi/n), k=0,\ldots,n$
be the Chebyshev points on the interval.
The Lagrange form of the 
interpolating polynomial is 
\begin{equation}
p(x) = 
\sum_{k=0}^n f(s_k)L_k(x),
\end{equation}
where $L_k$ are the Lagrange polynomials 
satisfying 
$L_k(s_{k^\prime})=\delta_{kk^\prime}$.
Among several expressions for $L_k$, 
we consider the barycentric 
form~\cite{berrut2004barycentric,
salzer1972lagrangian},
\begin{equation}
\label{eq:barycentric}
L_k(x) =
\frac{\displaystyle \frac{w_k}{x-s_k}}
{\displaystyle 
\sum_{k^\prime=0}^n
\frac{w_{k^\prime}}{x-s_{k^\prime}}},
\quad
w_k = (-1)^k
\begin{cases}
1/2, & k=0,n, \\ 1, & k=1:n-1.
\end{cases}
\end{equation}
Interpolation at the Chebyshev points
provides good uniform 
accuracy~\cite{boyd2001chebyshev,trefethen2019approximation},
and
the barycentric form permits 
stable and efficient
calculations~\cite{berrut2004barycentric,
higham2004numerical}.
This form is also scale-invariant in that
the same weights $w_k$
can be used for any interval $[a,b]$.

Two-dimensional
barycentric Lagrange interpolation has the form
\begin{equation}
p({\bf x}) = 
\sum_{\bf k}f({\bf s}_{\bf k})L_{\bf k}({\bf x}) =
\sum_{k_1=0}^n\sum_{k_2=0}^n 
f(s_{k_1}s_{k_2})L_{k_1}(x_1)L_{k_2}(x_2),
\end{equation}
where
${\bf x} = (x_1,x_2) \in [-1,1]^2$,
${\bf k} = (k_1,k_2)$ for $k_1,k_2 = 0:n$,
${\bf s}_{\bf k} = (s_{k_1},s_{k_2})$ 
are the tensor product Chebyshev points,
and
$L_{\bf k}({\bf x}) =
L_{k_1}(x_1)L_{k_2}(x_2)$.
This is implemented
on cubed sphere grid cells
using local angle 
coordinates~$(\xi,\eta)$~\cite{ronchi1996cubed},
where each cell on the sphere
is parameterized
by a square in~$(\xi,\eta)$,
and
the reference square $[-1,1]^2$ 
in the plane is mapped bilinearly 
to the cells on the sphere.
Henceforth for notational simplicity,
${\bf s}_{\bf k}$ 
shall denote the Chebyshev points 
in a cubed sphere grid cell $C$
(\cref{fig:eacs} shows a 
$5 \times 5$ example),
and as we shall see,
they serve as proxies for the 
target and source particles in $C$.
Next we explain how
barycentric Lagrange interpolation is used to
approximate interactions between
target and source particles on the sphere. 

\section{Interactions}
\label{sec:interactions}

\begin{figure}[htb]
    \centering
     \begin{subfigure}{0.45\textwidth}
        \includegraphics[width=\linewidth]{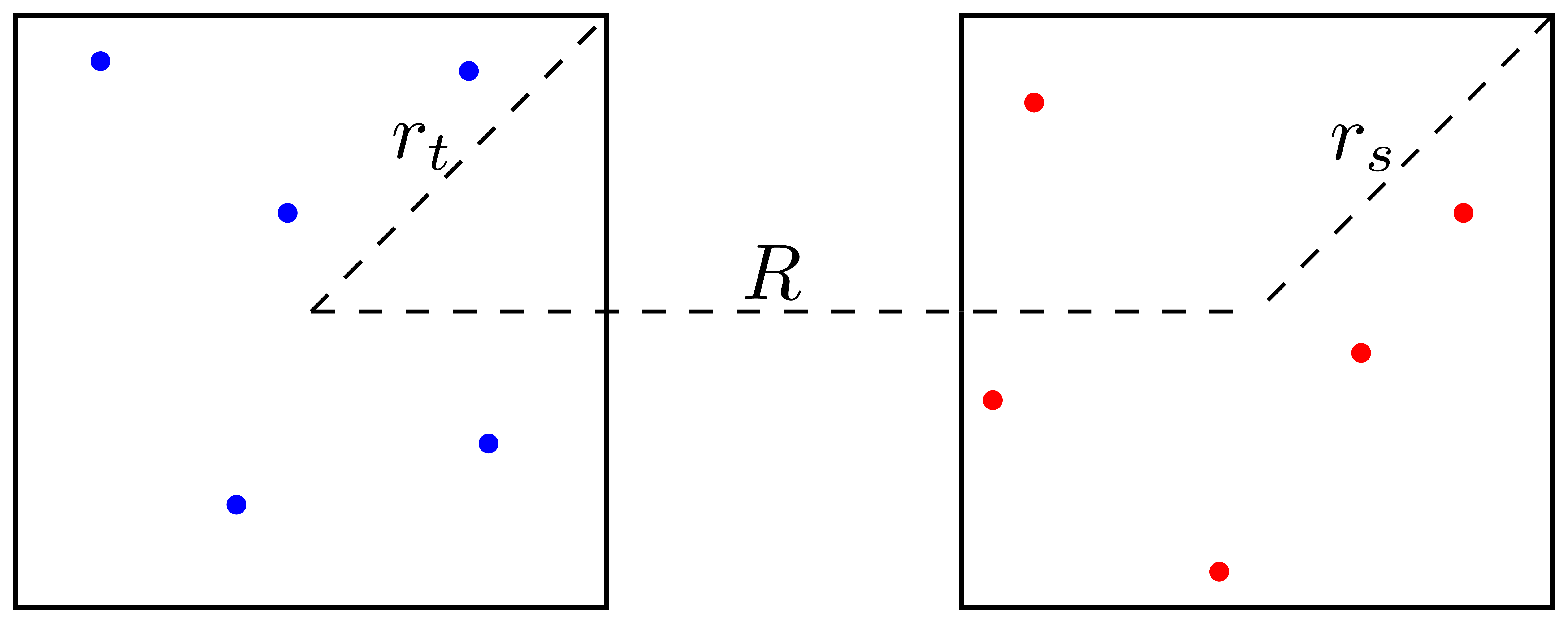}
        \caption{Particle-Particle (PP)}
        \label{fig:pp}
    \end{subfigure}
    \hspace{0.05\textwidth}
    \begin{subfigure}{0.45\textwidth}
        \includegraphics[width=\linewidth]{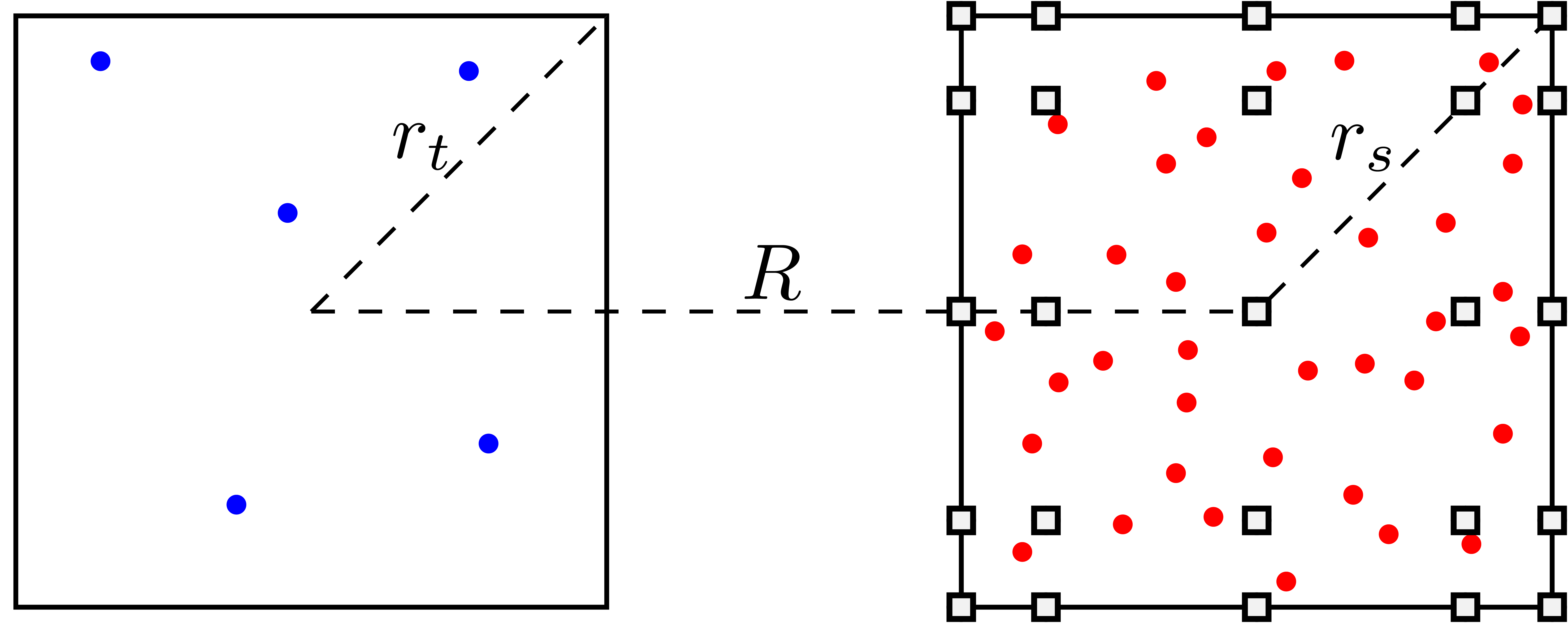}
        \caption{Particle-Cluster (PC)}
        \label{fig:pc}
    \end{subfigure}
    \begin{subfigure}{0.45\textwidth}
        \includegraphics[width=\linewidth]{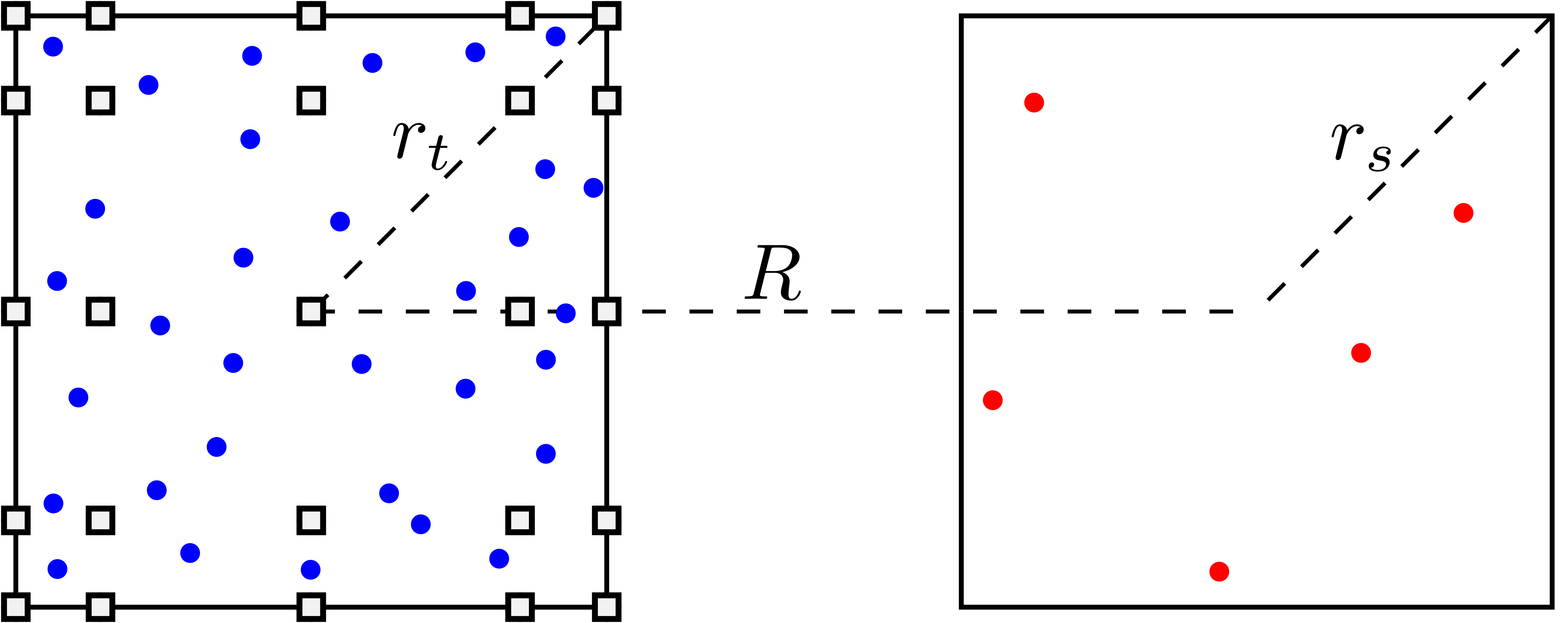}
        \caption{Cluster-Particle (CP)}
        \label{fig:cp}
    \end{subfigure}
    \hspace{0.05\textwidth}
    \begin{subfigure}{0.45\textwidth}
        \includegraphics[width=\linewidth]{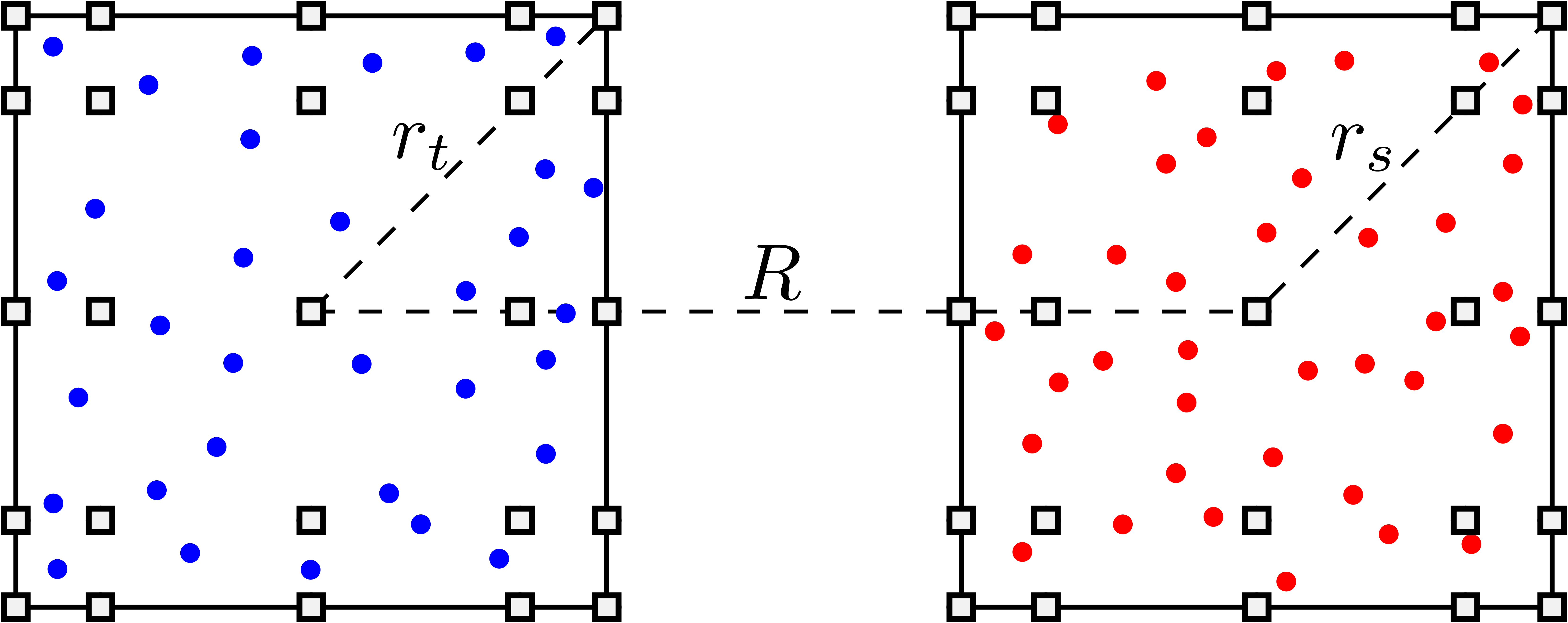}
        \caption{Cluster-Cluster (CC)}
        \label{fig:cc}
    \end{subfigure}
\vskip -15pt
\caption{Four types of interactions,
in each case the
target cluster $C_t$ is on the left 
(target particles in blue~\textcolor{blue}{$\bullet$})
and
the source cluster $C_s$ is on the right 
(source particles in red~\textcolor{red}{$\bullet$}), 
sample $5 \times 5$ grid of Chebyshev
proxy particles is shown (black~$\square$),
$R$ = distance between cluster centers, 
$r_t,r_s$ = radii of $C_t, C_s$.}
\label{fig:fourinteractions}
\end{figure}

\Cref{fig:fourinteractions}
depicts the four types of interactions used 
in the proposed fast summation methods,
referred to as
particle-particle (PP), 
particle-cluster (PC), 
cluster-particle (CP), 
and cluster-cluster (CC). 
The PP interactions are computed directly
with no approximation,
while the PC, CP, and CC interactions are
approximations obtained by
barycentric Lagrange interpolation 
on the cubed sphere grid cells.
In each case \Cref{fig:fourinteractions}
shows a target cluster $C_t$ with radius $r_t$
on the left
and
a source cluster $C_s$ with radius $r_s$
on the right,
as well as the distance $R$ between the
cluster centers.
Next we explain the interactions,
where ${\bf x}_i$ is a target particle in a
target cluster $C_t$
that interacts with the 
source particles ${\bf y}_j$ in a source cluster $C_s$.

\vskip 5pt
\noindent{\bf Particle-particle interaction.}
In a PP interaction (\Cref{fig:pp}),
$C_t$ and $C_s$ have few particles. 
In this case the potential at a 
target particle ${\bf x}_i$ in $C_t$
due to the source particles ${\bf y}_j$ in $C_s$ is
computed directly with no approximation,
\begin{equation}\label{eq:PP}
\phi(\mathbf{x}_i,C_t,C_s) =
\sum_{\mathbf{y}_j \in C_s}K(\mathbf{x}_i,\mathbf{y}_j)w_j,
\end{equation} 
where $w_j$ are the source particle weights.
The expression~\cref{eq:PP}
is called a PP potential.

\vskip 5pt
\noindent{\bf Particle-cluster interaction.}
In a PC interaction (\cref{fig:pc}), 
$C_t$ has few particles 
and $C_s$ has many particles. 
In this case the potential at a 
target particle ${\bf x}_i$ in $C_t$
due to the source particles ${\bf y}_j$ in $C_s$ is
approximated by interpolating the kernel
in the source variable,
\begin{equation}
\phi(\mathbf{x}_i,C_t,C_s) =
\sum_{\mathbf{y}_j \in C_s}
K(\mathbf{x}_i,\mathbf{y}_j)w_j \approx
\sum_{\mathbf{y}_j \in C_s}
\sum_{\bf k}
K(\mathbf{x}_i,\mathbf{s}_{\bf k})
L_{\bf k}({\bf y}_j)w_j,
\end{equation}
where
${\bf k} = (k_1,k_2)$ for $k_1,k_2 = 0:n$,
${\bf s}_{\bf k}$ are the 
Chebyshev points in $C_s$
(also called proxy source particles),
and
$L_{\bf k}({\bf y}) =
L_{k_1}(y_1)L_{k_2}(y_2)$ for ${\bf y} = (y_1,y_2)$.
This leads to the PC potential,
\begin{equation}
\label{eq:PC}
\phi({\bf x}_i,C_t,\overline{C_s}) =
\sum_{\bf k}
K(\mathbf{x}_i,\mathbf{s}_{\bf k})
\overline{w}_{\bf k}, \quad
\overline{w}_{\bf k} = \sum_{\mathbf{y}_j \in C_s}
L_{\bf k}({\bf y}_j)w_j,
\end{equation} 
where 
$\overline{C_s}$ indicates that 
interpolation occurred in the source variable,
and $\overline{w}_{\bf k}$ are the 
proxy weights of $C_s$.
The proxy weights $\overline{w}_{\bf k}$ 
are independent of $\mathbf{x}_i$, 
and once computed,
they are stored and reused for 
PC interactions between other 
target particles and~$C_s$.
Instead of computing $\overline{w}_{\mathbf{k}}$ 
directly for each source cluster $C_s$, 
a more efficient upward pass described below 
is used.

\vskip 5pt
\noindent{\bf Cluster-particle interaction.}
In a CP interaction (\cref{fig:cp}), 
$C_t$ has many particles 
and $C_s$ has few particles. 
In this case the potential at a 
target particle ${\bf x}_i$ in $C_t$
due to the source particles ${\bf y}_j$ in $C_s$ 
is approximated by interpolating the 
kernel in the target variable,
\begin{subequations}
\begin{align}
\phi(\mathbf{x}_i,C_t,C_s) =
\sum_{\mathbf{y}_j \in C_s}
K(\mathbf{x}_i,\mathbf{y}_j)w_j
&\approx \sum_{\mathbf{y}_j \in C_s}
\sum_{\bf m}
K(\mathbf{t}_{\bf m},\mathbf{y}_j)
L_{\bf m}({\bf x}_i)w_j \\
\label{eq:CP_interpolation}
&= \sum_{\bf m}\phi({\bf t}_{\bf m},\overline{C_t},C_s)
L_{\bf m}({\bf x}_i),
\end{align}
\end{subequations}
where
${\bf t}_{\bf m}$ are the 
Chebyshev points in $C_t$
(also called proxy target particles),
${\bf m} = (m_1,m_2)$ for $m_1,m_2 = 0:n$,
$L_{\bf m}({\bf y}) =
L_{m_1}(y_1)L_{m_2}(y_2)$ for ${\bf y} = (y_1,y_2)$,
and
the CP proxy potentials are defined by
\begin{equation}
\label{eq:CP}
\phi({\bf t}_{\bf m},\overline{C_t},C_s) =
\sum_{\mathbf{y}_j \in C_s}
K(\mathbf{t}_{\bf m},\mathbf{y}_j)w_j.
\end{equation}
The notation $\overline{C}_t$ 
indicates that interpolation occurred 
in the target variable. 
Instead of interpolating the
CP proxy potentials directly from $\mathbf{t}_{\mathbf{m}}$ to $\mathbf{x}_i$
by~\cref{eq:CP_interpolation},  
a more efficient downward pass described below
is used.

\vskip 5pt
\noindent{\bf Cluster-cluster interaction.}
In a CC interaction (\cref{fig:cc}), 
$C_t$ and $C_s$ have many particles. 
In this case the potential at a 
target particle ${\bf x}_i$ in $C_t$
due to the source particles ${\bf y}_j$ in $C_s$ is
approximated by interpolating the kernel
in the target and source variables,
\begin{subequations}
\begin{align}
\phi(\mathbf{x}_i,C_t,C_s) =
\sum_{\mathbf{y}_j \in C_s}
K(\mathbf{x}_i,\mathbf{y}_j)w_j
&\approx
\sum_{\mathbf{y}_j \in C_s}
\sum_{\bf m}\sum_{\bf k}
K(\mathbf{t}_{\bf m},\mathbf{s}_{\bf k})
L_{\bf m}({\bf x}_i)L_{\bf k}({\bf y}_j)w_j \\
\label{eq:CC_interpolation}
&= \sum_{\bf m}
\phi({\bf t}_{\bf m},
\overline{C_t},\overline{C_s})
L_{\bf m}({\bf x}_i),
\end{align}
\end{subequations}
where the CC proxy potentials are
\begin{equation}
\label{eq:CC}
\phi({\bf t}_{\bf m},
\overline{C_t},\overline{C_s}) =
\sum_{\bf k}
K({\bf t}_{\bf m},{\bf s}_{\bf k})
\overline{w}_{\bf k},
\quad
\overline{w}_{\bf k} = \sum_{\mathbf{y}_j \in C_s}
L_{\bf k}({\bf y}_j)w_j.
\end{equation} 
The notation
$\overline{C_t}, \overline{C_s}$ 
indicates that interpolation occurred
in the target and source variables,
and
the proxy weights $\overline{w}_{\bf k}$
are the same as in~\eqref{eq:PC}.
Instead of interpolating the
CC proxy potentials directly from $\mathbf{t}_{\mathbf{m}}$ to $\mathbf{x}_i$
by~\cref{eq:CC_interpolation},  
a more efficient downward pass described below
is used.

This concludes the description of the
four types of interactions.
Note that computing the 
PP and PC potentials in
\eqref{eq:PP} and \eqref{eq:PC}, 
and the
CP and CC proxy potentials in 
\eqref{eq:CP} and \eqref{eq:CC}
call for direct summation
of interactions between
target particles ${\bf x}_i,{\bf t}_{\bf m}$
and 
source particles ${\bf y}_j,{\bf s}_{\bf k}$
with weights $w_j,\overline{w}_{\bf k}$;
this is important because
direct summation can be efficiently parallelized.

\section{Upward pass}
\label{sec:upward}

The upward pass computes the proxy weights $\overline{w}_{\bf k}$ defined in
\eqref{eq:PC} for each source cluster $C_s$.
Given a parent cluster $C_s^p$ 
with child clusters $C_s^c, c = 1:4$,
the proxy weight of the parent is
\begin{equation}\label{eq:upwardpass1}
\overline{w}_{\bf k}(C_s^p) =
\sum_{\mathbf{y}_j\in C_s^p}L_{\bf k}^p(\mathbf{y}_j)w_j =
\sum_{c=1}^4
\sum_{\mathbf{y}_j\in C_s^c}
L_{\bf k}^p(\mathbf{y}_j)w_j,
\end{equation}
where $L_{\bf k}^p$ are the Lagrange polynomials
of the parent.
Next let $\mathbf{s}_{\bf m}^c$ 
denote the Chebyshev points 
in the child cluster $C_s^c$
with Lagrange polynomials $L_{\bf m}^c$.
Since polynomial interpolation 
is exact for polynomials, it follows that
\begin{equation}
\label{eq:interpolation_identity}
L_{\bf k}^p(\mathbf{y}_j) =
\sum_{\bf m}
L_{\bf k}^p(\mathbf{s}_{\bf m}^c)
L_{\bf m}^c(\mathbf{y}_j),
\quad c = 1:4,
\end{equation}
and inserting this in~\eqref{eq:upwardpass1}
yields
\begin{equation}
\label{eq:upwardpass2}
\overline{w}_{\bf k}(C_s^p) =
\sum_{c=1}^4\sum_{\bf m}
L_{\bf k}^p(\mathbf{s}_{\bf m}^c)
\overline{w}_{\bf m}(C_s^c), \quad
\overline{w}_{\bf m}(C_s^c) =
\sum_{\mathbf{y}_j\in C_s^c}
L_{\bf m}^c(\mathbf{y}_j)w_j.
\end{equation}
Hence the parent proxy weights
$\overline{w}_{\bf k}(C_s^p)$ 
can be obtained from the child proxy weights 
$\overline{w}_{\bf m}(C_s^c)$. 
The proxy weights of the leaf clusters
are computed directly using~\eqref{eq:PC},
and
they are combined to compute the
proxy weights of the parent clusters
using~\eqref{eq:upwardpass2},
and
so on for higher levels in the tree. 
The proxy weights are 
computed once for each source cluster,
and
then stored and used as needed.
The upward pass described here is
the same as in the
FMM~\cite{cheng1999fast,fmmgreengard},
however adapted to polynomial 
interpolation~\cite{bldtt}.

\section{Cubed Sphere Tree Code}
\label{sec:cstc}

The CSTC follows the
Barnes-Hut approach~\cite{barneshut}, modified to use
barycentric Lagrange interpolation 
on cubed sphere grid cells
rather than monopole approximations
on rectangular boxes in $\mathbb{R}^3$.
The potential $\phi({\bf x}_i)$
has contributions from 
near-field PP interactions 
and
far-field PC interactions.

\Cref{alg:cstc} presents the CSTC pseudocode.
Lines~1 and~2 specify
the input particle data and numerical parameters.
Line~3 initializes the potentials 
$\phi({\bf x}_i)$.
Line~4 builds the source tree
from the source particles $\mathbf{y}_j$,
where the source panels $C_s$ 
are cubed sphere grid cells at different
levels of refinement.
Line~5 performs the upward pass to compute the 
proxy source weights $\overline{w}_{\bf k}$.
Line~6 starts a loop over the 
target particles ${\bf x}_i$
and
Line~7 starts a downward traversal
of the source tree through
source clusters $C_s$.
Line~8 checks whether ${\bf x}_i$ and $C_s$ 
are well-separated;
the criterion is $r/R < {\rm MAC}$,
where 
$r$ is the radius of $C_s$, 
$R$ is the distance between 
$\mathbf{x}_i$ and the center of $C_s$,
and
MAC is a user-specified parameter
(multipole acceptance criterion~\cite{barneshut}).
If $\mathbf{x}_i$ and $C_s$ are well-separated,
the size of $C_s$ is checked;
if $|C_s| > N_0$,
the PC~approximation~\eqref{eq:PC} 
is computed (Line~10),
and
if $|C_s| \le N_0$,
the direct PP~interaction~\eqref{eq:PP} 
is computed (Line 12). 
If $\mathbf{x}_i$ and $C_s$ are not well-separated, 
the algorithm considers interactions
between $\mathbf{x}_i$ and the child clusters
$C_s^c, c=1:4$ (Line 14).
Line~15 outputs the potentials $\phi({\bf x}_i)$.

\begin{algorithm}[htb]
\caption{Cubed Sphere Tree code}
\label{alg:cstc}
\begin{algorithmic}[1]
\STATE{{\bf input}: 
target particles ${\bf x}_i$,
source particles ${\bf y}_j$,
source weights $w_j$}
\STATE{{\bf input}:
interpolation degree $n$,
MAC parameter,
maximum leaf size $N_0$}
\STATE{initialize 
potentials $\phi({\bf x}_i)=0$}
\STATE{build source tree from source
particles ${\bf y}_j$}
\STATE{upward pass to compute proxy source weights $\overline{w}_{\bf k}$}
\FOR{target particle ${\bf x}_i, i=1:N$}
\FOR{cluster $C_s$ in source tree}
\IF{$\mathbf{x}_i$ and $C_s$ are well-separated}
\IF{$|C_s| > N_0$} 
\STATE{compute PC approximation~\eqref{eq:PC} between ${\bf x}_i$ and $C_s$}
\ELSE
\STATE{compute direct PP~interaction~\eqref{eq:PP}
between ${\bf x}_i$ and $C_s$}
\ENDIF
\ELSE
\STATE{return to Line~7 and interact with
child clusters $C_s^c, c=1:4$}
\ENDIF
\ENDFOR
\ENDFOR
\STATE{{\bf output}: 
potentials $\phi({\bf x}_i)$}
\end{algorithmic}
\end{algorithm}

\section{Cubed Sphere Fast Multipole Method}
\label{sec:csfmm}

The CSFMM uses upward and downward passes 
as in the FMM~\cite{cheng1999fast,fmmgreengard}
with two modifications~\cite{bldtt},
(1) the kernel approximations use
barycentric Lagrange interpolation
on cubed sphere grid cells
rather than analytic multipole expansions
on boxes in $\mathbb{R}^3$,
(2) the interactions are determined by
dual tree traversal~\cite{appel1985efficient}
rather than FMM interaction lists. 
In the CSFMM the potential $\phi({\bf x}_i)$
has contributions from 
near-field PP interactions, 
and
far-field PC, CP and CC interactions;
the PP and PC interactions contribute 
potentials at target particles ${\bf x}_i$,
while the
CP and CC interactions 
compute proxy potentials 
at proxy target particles
${\bf t}_{\bf m}$
which are 
accumulated in a downward pass 
through the target tree and 
interpolated to
potentials at target particles ${\bf x}_i$
in the leaves.

\Cref{alg:csfmm} 
presents the CSFMM pseudocode.
Lines~1 and~2 specify
the input particle data and numerical parameters.
Line~3 initializes the potentials
$\phi({\bf x}_i)$.
Line~4 builds two trees, 
a target tree with clusters~$C_t$
for the target particles $\mathbf{x}_i$ 
and 
a source tree with clusters $C_s$
for the source particles $\mathbf{y}_j$.
Line~5 performs the upward pass to compute
the proxy source weights $\overline{w}_{\bf k}$.
Line~6 starts the dual tree traversal
sketched in Lines~7-18.
For a given target cluster $C_t$ 
and source cluster $C_s$,
if $C_t$ and $C_s$ are well-separated
and 
have sufficiently many particles, 
then a PC, CP or CC interaction is performed
depending on the sizes of $C_t$ and $C_s$.
Line~15 accumulates the 
CP and CC proxy potentials computed
in Line 12 and Line~14.
If $C_t$ and $C_s$ are not well-separated,
the code returns to Line~7
and
considers interactions with the
children of $C_t$ or $C_s$.
Line~19 adds the PP interactions to the potentials.
Line~20 performs the downward pass 
to interpolate the 
accumulated proxy potentials
from the proxy target particles~${\bf t}_{\bf m}$
to the target particles ${\bf x}_i$.
Line~21 outputs the potentials $\phi({\bf x}_i)$.

\begin{algorithm}[htb]
\caption{Cubed Sphere Fast Multipole Method}
\label{alg:csfmm}
\begin{algorithmic}[1]
\STATE{{\bf input}: target particles ${\bf x}_i$,
source particles ${\bf y}_j$, 
source weights $w_j$}
\STATE{{\bf input}:
interpolation degree $n$,
MAC parameter,
maximum leaf size $N_0$}
\STATE{initialize potentials 
$\phi({\bf x}_i)=0$}
\STATE{build target tree with clusters $C_t$
and source tree with clusters $C_s$}
\STATE{upward pass to compute 
proxy source weights $\overline{w}_k$}
\STATE{start dual tree traversal}
\FOR{target cluster $C_t$ and source cluster $C_s$}
\IF{$C_t$ and $C_s$ are well-separated
and $|C_t| > N_0$ or $|C_s| > N_0$}
\IF{$|C_t|<N_0$ and $|C_s|>N_0$} 
\STATE{add PC interaction
$\phi({\bf x}_i,C_t,\overline{C_s})$
in~\eqref{eq:PC} to
$\phi({\bf x}_i)$ for ${\bf x}_i \in C_t$}
\ENDIF
\IF{$|C_t|>N_0$ and $|C_s|<N_0$} 
\STATE{compute CP proxy potentials
$\phi({\bf t}_{\bf m},\overline{C_t},C_s)$
in~\eqref{eq:CP}
for ${\bf t}_{\bf m} \in C_t$}
\ELSE
\STATE{compute CC proxy potentials
$\phi({\bf t}_{\bf m},\overline{C_t},\overline{C_s})$
in~\eqref{eq:CC}
for ${\bf t}_{\bf m} \in C_t$}
\ENDIF
\STATE{accumulate 
CP and CC proxy potentials in
$\phi(\mathbf{t}_{\bf m},C_t)$ 
for ${\bf t}_{\bf m} \in C_t$}
\ENDIF
\IF{$C_t$ and $C_s$ are not well-separated
and $|C_t|>N_0$ or $|C_s|>N_0$}
\STATE{return to Line~7 and interact with the
children of the larger of $C_t$ and $C_s$}
\ELSE
\STATE{add PP interactions~\eqref{eq:PP} 
to potentials $\phi({\bf x}_i)$
for ${\bf x}_i \in C_t$}
\ENDIF
\ENDFOR
\STATE{downward pass to
interpolate accumulated proxy potentials
from ${\bf t}_{\bf m}$ to $\mathbf{x}_i$}
\STATE{{\bf output}: 
potentials $\phi({\bf x}_i)$}
\end{algorithmic}
\end{algorithm}

Some details need explanation.
The dual tree traversal in Line~6
starts by cycling over pairs of
root clusters from the target and source trees.
In Line~8,
a target cluster $C_t$
and
source cluster $C_s$ 
are considered to be well-separated if
$(r_t+r_s)/R<\mathrm{MAC}$, 
where 
$r_t$,$r_s$ are the radii of $C_t$, $C_s$, 
$R$ is the distance between the centers of $C_t$ and $C_s$, 
and 
$\mathrm{MAC}$ is a user-specified parameter.
Line 15 accumulates the 
CP and CC proxy potentials,
\begin{equation}  
\label{eq:CP_CC_accumulated}
\phi(\mathbf{t}_{\bf m},C_t) =
\sum_{\mathrm{CP}}
\phi({\bf t}_{\bf m},\overline{C_t},C_s) +
\sum_{\mathrm{CC}}
\phi({\bf t}_{\bf m},
\overline{C_t},\overline{C_s}),
\quad {\bf t}_{\bf m} \in C_t,
\end{equation}
where the sums are taken over the
source clusters $C_s$
that interacted with the target cluster $C_t$
in Line~12 or Line~14.
After the dual tree traversal finishes (Line~19),
the potentials $\phi({\bf x}_i)$
contain contributions from all the
PP and PC interactions.
The accumulated proxy potentials
$\phi({\bf t}_{\bf m},C_t)$
from CP and CC interactions
could be interpolated directly
to the target particles ${\bf x}_i \in C_t$,
but this is accomplished
more efficiently by the downward pass 
(Line 20),
which is described in~\Cref{sec:downward}.

\subsection{Downward pass}
\label{sec:downward}

For simplicity assume the tree has two levels,
a parent level and a child level.
Consider a target particle 
$\mathbf{x}_i \in C_t^c \subset C_t^p$, 
where
$C_t^c$ is a child cluster
of the parent cluster $C_t^p$.
Then the potential at $\mathbf{x}_i$ 
due to the accumulated
CP and CC interactions~\cref{eq:CP_CC_accumulated}
of the child and parent can be written as
\begin{equation} 
\label{eq:downward_1}
\phi(\mathbf{x}_i,C_t^c+C_t^p) =
\sum_{\mathbf{n}}
\phi({\bf t}_{\bf n}^c,C_t^c)
L_{\bf n}^c(\mathbf{x}_i) +
\sum_{\bf m}
\phi({\bf t}_{\bf m}^p,C_t^p)
L_{\bf m}^p(\mathbf{x}_i),
\end{equation}
where 
the 1st~sum interpolates from 
child proxy particles ${\bf t}_{\bf n}^c$ to ${\bf x}_i$
and
the 2nd~sum interpolates from 
parent proxy particles ${\bf t}_{\bf m}^p$ to ${\bf x}_i$.
Alternatively, writing the
identity~\eqref{eq:interpolation_identity}
in the form
\begin{equation}
L_{\bf m}^p(\mathbf{x}_i) =
\sum_{\bf n}
L_{\bf m}^p({\bf t}_{\bf n}^c)
L_{\bf n}^c(\mathbf{x}_i),
\end{equation}
and
inserting this in~\cref{eq:downward_1} yields
\begin{equation}   
\label{eq:downward_2}
\phi(\mathbf{x}_i,C_t^c+C_t^p) =
\sum_{\mathbf{n}}
\bigg(
\phi(\mathbf{t}_{\mathbf{n}}^c,C_t^c) +
\sum_{\mathbf{m}}
\phi(\mathbf{t}_{\mathbf{m}}^p,C_t^p)
L_{\mathbf{m}}^p(\mathbf{t}_{\mathbf{n}}^c)
\bigg)
L_{\mathbf{n}}^c(\mathbf{x}_i),
\end{equation}
which replaces the interpolation from
${\bf t}_{\bf m}^p$ to ${\bf x}_i$
in the 2nd~sum of~\cref{eq:downward_1}
by interpolation from
${\bf t}_{\bf m}^p$ to ${\bf t}_{\bf n}^c$
and then from 
${\bf t}_{\rm n}^c$ to ${\bf x}_i$.
This idea generalizes to a tree with
more than two levels,
where parent proxy potentials are interpolated
to child proxy potentials
through all levels in the tree
before finally being interpolated to
potentials at the leaf level.
The downward pass described here is
the same as in the
FMM~\cite{cheng1999fast,fmmgreengard},
however adapted to polynomial 
interpolation~\cite{bldtt}.

\subsection{Related methods}

Here we briefly discuss 
the key features that distinguish the CSFMM
from prior related work.
The CSFMM is most closely related to the 
recently developed BLDTT~\cite{bldtt}.
Both methods approximate the kernel using
barycentric Lagrange interpolation,
but BLDTT applies to 
kernels defined in Euclidean space,
whereas CSFMM applies to
kernels defined on the sphere.
Even if the target and source points 
lie on the sphere,
the interpolation points utilized in BLDTT
lie off the sphere.
Hence BLDTT cannot be used
if the kernel is defined only on the sphere.
An example is the 
approximate SAL kernel~\cref{eq:G_SAL_approx}
which requires 
$\sqrt{2(1 - {\bf x} \cdot {\bf y})}$;
this is well-defined for 
points ${\bf x}, {\bf y}$ on the sphere where 
${\bf x} \cdot {\bf y} = \cos\varphi \le 1$,
but BLDTT would require 
kernel evaluations at points off the sphere
where in some cases ${\bf x} \cdot {\bf y} > 1$.
The same issue will arise for other methods
if they employ kernel evaluations at proxy points off the 
sphere~\cite{fong2009black,
xing2020interpolative,ying2004kernel},
however
CSFMM avoids the problem by requiring
kernel evaluations only at points on the sphere.

Among numerical methods
for the Poisson equation on a sphere,
one example is a spherical Fast Multipole Method (sFMM) 
for gravitational lensing simulations
that uses an analytic series expansion for the
spherical Laplace Green's 
function~\cref{eq:G_L}~\cite{suo2023spherical}.
By contrast,
CSFMM does not employ analytic series expansions
and
requires only the ability
to evaluate the kernel at general points on the sphere.
In that sense CSFMM is kernel-independent
and can be broadly applied
within the class of spherical kernels
as shown below for the kernels discussed 
in \Cref{sec:sphereconv}.

Finally we mention a spherical tree code (STC)
in which the kernel was approximated using
polynomial interpolation on 
spherical triangles~\cite{chen2024fast}.
The polynomial coefficients were obtained by 
solving a linear system,
but the system becomes ill-conditioned 
as the interpolation degree increases
and
this limits the accuracy that can be attained.
The barycentric Lagrange approach in CSFMM
enables accurate stable interpolation 
without having to 
solve linear systems~\cite{berrut2004barycentric,
higham2004numerical}.

\section{Implementation details}
\label{sec:details}

The CSFMM and CSTC were applied to the problems 
described in~\Cref{sec:sphereconv}. 
The code was written in double precision C++
and
was compiled using the Intel OneAPI compiler
with optimization level O2.
The calculations were done on the 
NCAR Derecho system~\cite{derecho},
where each node has 128 cores 
corresponding to two 3rd~generation AMD EPYC 7763 Milan processors.
\Cref{tbl:grids} 
gives the parameters for computing
the $N$-body sum~\cref{eq:nbodysum}
with the three
spherical partitions in~\Cref{fig:grids}.
The angular grid spacing for the icosahedral and cubed sphere grids is approximated as $\Delta\varphi = \sqrt{4\pi/N} \cdot 180/\pi$
with units in degrees. 
At each successive level of refinement
the particle count $N$ increases by a factor
of approximately four
and the
angular grid spacing $\Delta\varphi$
decreases by one half.
The code for these tests is available online~\cite{CSFMM}. 

\begin{table}[htb]
\centering
\begin{tabular}{|c|rc|rc|rc|}
\hline
 &
\multicolumn{2}{|c|}{icosahedral} &
\multicolumn{2}{c|}{cubed sphere} &
\multicolumn{2}{c|}{\begin{tabular}[c]{@{}c@{}} latitude-longitude \\ \end{tabular}} \\ 
\hline
$L$ &
\multicolumn{1}{|c|}{$N$} &
\multicolumn{1}{|c|}{$\Delta\varphi$} &
\multicolumn{1}{|c|}{$N$} &
\multicolumn{1}{|c|}{$\Delta\varphi$} &
\multicolumn{1}{|c|}{$N$} &
\multicolumn{1}{|c|}{$\Delta\varphi$} \\
\hline
\multicolumn{1}{|c|}{4} &
\multicolumn{1}{r|}{2562} & $4.01^{\circ}$ &
\multicolumn{1}{r|}{1536} & $5.18^{\circ}$ &
\multicolumn{1}{r|}{4050} & $4^{\circ}$ \\ 
\hline
\multicolumn{1}{|c|}{5} &
\multicolumn{1}{r|}{10242} & $2.01^{\circ}$ &
\multicolumn{1}{r|}{6144} & $2.59^{\circ}$ &
\multicolumn{1}{r|}{16200} & $2^{\circ}$ \\ 
\hline
\multicolumn{1}{|c|}{6} &
\multicolumn{1}{r|}{40962} & $1.00^{\circ}$ &
\multicolumn{1}{r|}{24576} & $1.30^{\circ}$ &
\multicolumn{1}{r|}{64800} & $1^{\circ}$ \\ 
\hline
\multicolumn{1}{|c|}{7} &
\multicolumn{1}{r|}{163842} & $0.502^{\circ}$ &
\multicolumn{1}{r|}{98304} & $0.648^{\circ}$ &
\multicolumn{1}{r|}{259200} & $0.5^{\circ}$ \\ 
\hline
\multicolumn{1}{|c|}{8} &
\multicolumn{1}{r|}{655362} & $0.251^{\circ}$ &
\multicolumn{1}{r|}{393216} & $0.324^{\circ}$ &
\multicolumn{1}{r|}{1036800} & $0.25^{\circ}$ \\ 
\hline
\multicolumn{1}{|c|}{9} &
\multicolumn{1}{r|}{2621442} & $0.125^{\circ}$ &
\multicolumn{1}{r|}{1572864} & $0.162^{\circ}$ &
\multicolumn{1}{r|}{4147200} & $0.125^{\circ}$ \\ 
\hline
\end{tabular}
\caption{
Parameters for computing $N$-body 
sum~\cref{eq:nbodysum} using the three
spherical partitions in~\Cref{fig:grids}
showing grid level $L$,
particle count $N$,
approximate angular grid spacing $\Delta\varphi$.}
\label{tbl:grids}
\end{table}

\section{Numerical results: error}
\label{sec:error}

In discussing the error
we consider three solutions;
$\phi_{\rm EX}$ is the reference exact solution,
and
$\phi_{\rm DS}, \phi_{\rm FS}$ 
use the midpoint rule to compute the integral
transform~\cref{eq:integral_transform},
where 
$\phi_{\rm DS}$ uses direct summation
to compute the $N$-body sum~\cref{eq:nbodysum}, 
and $\phi_{\rm FS}$ uses a
fast sum method
(CSFMM by default or CSTC if indicated).
The discretization 
error $|\phi_{\rm DS}-\phi_{\rm EX}|$
is controlled by the particle count $N$,
while the fast sum approximation error 
$|\phi_{\rm FS}-\phi_{\rm DS}|$ is controlled by the
MAC parameter
and
interpolation degree~$n$;
the maximum leaf size is set to $N_0=4n^2$,
which is roughly the number of proxy particles
in the four children of a parent cluster.
The fast sum approximation error can be reduced by
decreasing the MAC parameter
and
increasing the degree $n$,
but this raises the computational cost; 
we shall examine the effect of varying these parameters,
but unless stated otherwise,
the CSFMM and CSTC use 
$\mathrm{MAC}=0.7$ and degree $n=6$.
We also examine the fast sum solution
error $|\phi_{\rm FS}-\phi_{\rm EX}|$,
which includes
discretization error and
fast sum approximation error.

\subsection{Spatial distribution of error}

\begin{figure}[htb]
\includegraphics[width=\textwidth]{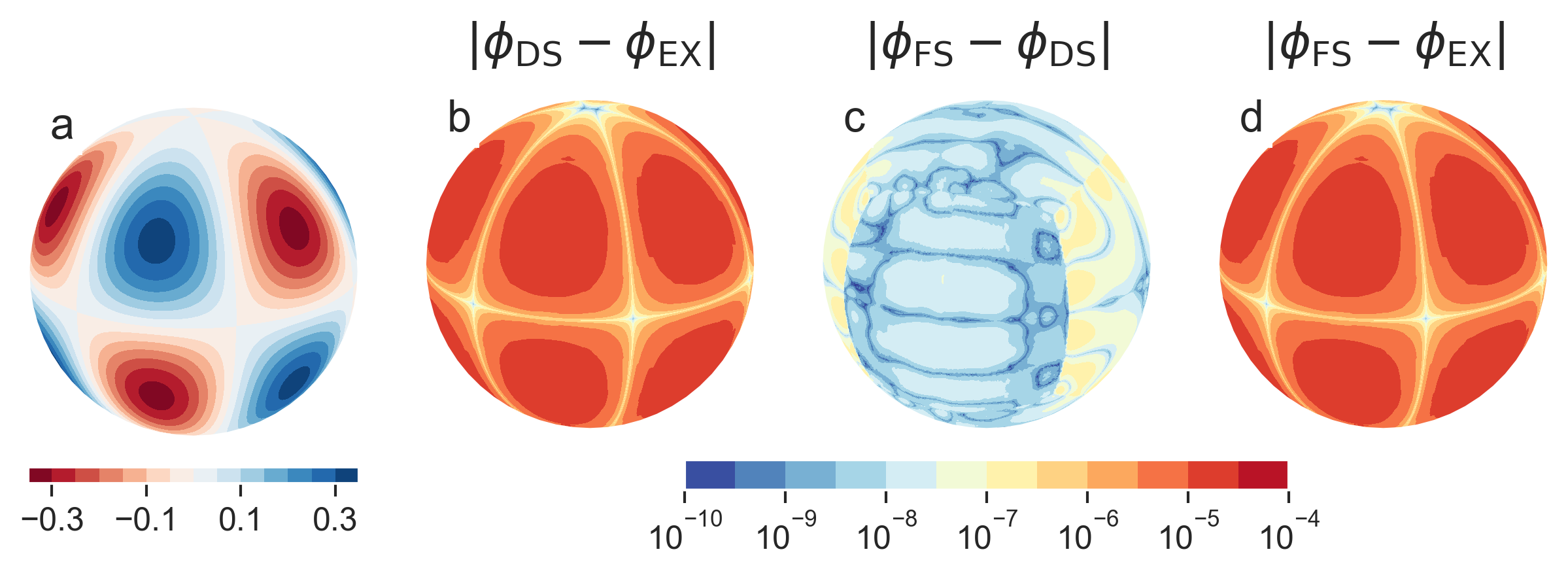}
\includegraphics[width=\textwidth]{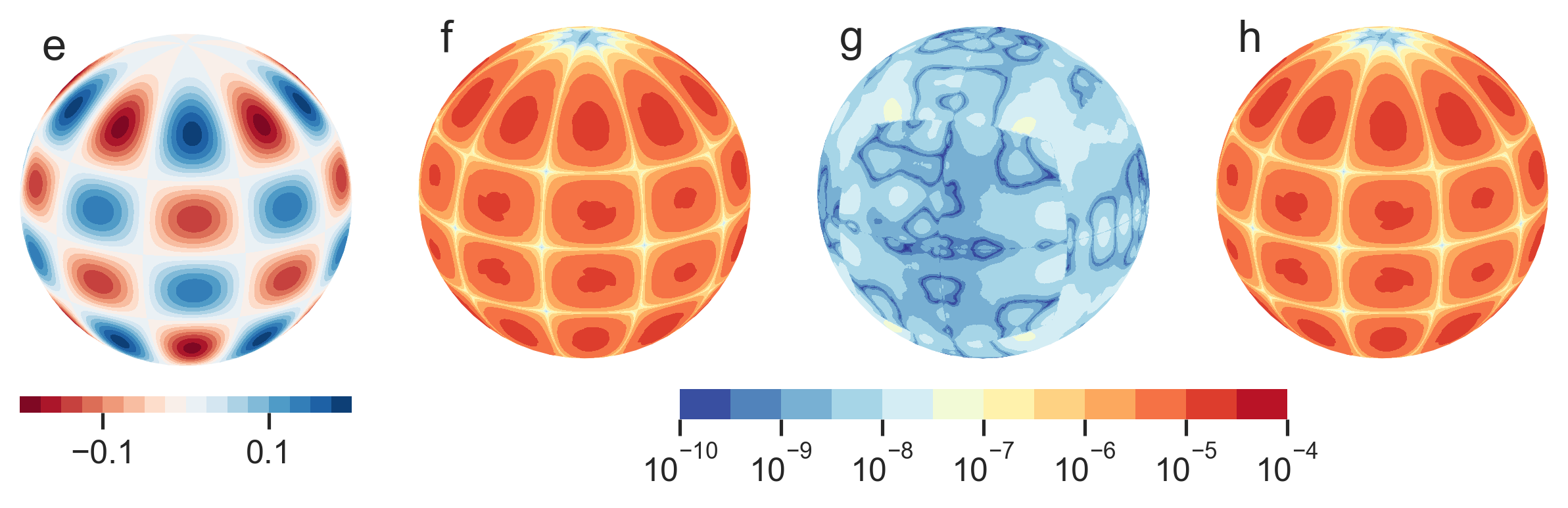}
\vskip -10pt
\caption{
Poisson equation~\cref{eq:Poisson_equation}, 
the integral in \cref{eq:G_L_1}
is computed for
real spherical harmonics
(a) $Y_{4,3}$,
(e) $Y_{8,5}$, 
icosahedral grid with $N=163842$ points,
(b,f) discretization error 
$|\phi_{\rm DS}-\phi_{\rm EX}|$,
(c,g) CSFMM approximation error 
$|\phi_{\rm FS}-\phi_{\rm DS}|$,
(d,h) CSFMM solution error 
$|\phi_{\rm FS}-\phi_{\rm EX})$,
CSFMM with MAC = 0.7, degree $n=6$.}
\label{fig:poissonspatial}
\end{figure}

This test considers the 
Poisson equation~\cref{eq:Poisson_equation}
in two cases where the exact solution is a
real spherical harmonic, 
$Y_{4,3}$ 
in~\cref{fig:poissonspatial}a
and
$Y_{8,5}$ 
in~\cref{fig:poissonspatial}e.
The numerical solution is obtained by
discretizing the integral
in~\cref{eq:G_L_1}
on the icosahedral grid with $N=163842$ points.
Figures~\ref{fig:poissonspatial}b,
\ref{fig:poissonspatial}f
show the spatial distribution of the 
discretization error
$|\phi_{\rm DS}-\phi_{\rm EX}|$, 
which is correlated with the exact solution,
and is largest at the
local maxima and minima of the exact solution,
as expected for the midpoint rule.
Figures~\ref{fig:poissonspatial}c,
\ref{fig:poissonspatial}g
show the CSFMM approximation error 
$|\phi_{\rm FS}-\phi_{\rm DS}|$,
which is 
several orders of magnitude smaller than the
discretization error,
with some minor imprint from the 
cubed sphere grid used in the CSFMM interpolation,
but uncorrelated with the exact solution. 
Figures~\ref{fig:poissonspatial}d,
\ref{fig:poissonspatial}h
show the CSFMM solution error
$|\phi_{\rm FS}-\phi_{\rm EX}|$,
confirming that the CSFMM approximation error
is negligible compared to the
discretization error.
Similar results were obtained using the 
cubed sphere and latitude-longitude grids.

\subsection{Error versus particle count}

\begin{figure}[h!]
\centering
\includegraphics[width=0.80\linewidth]{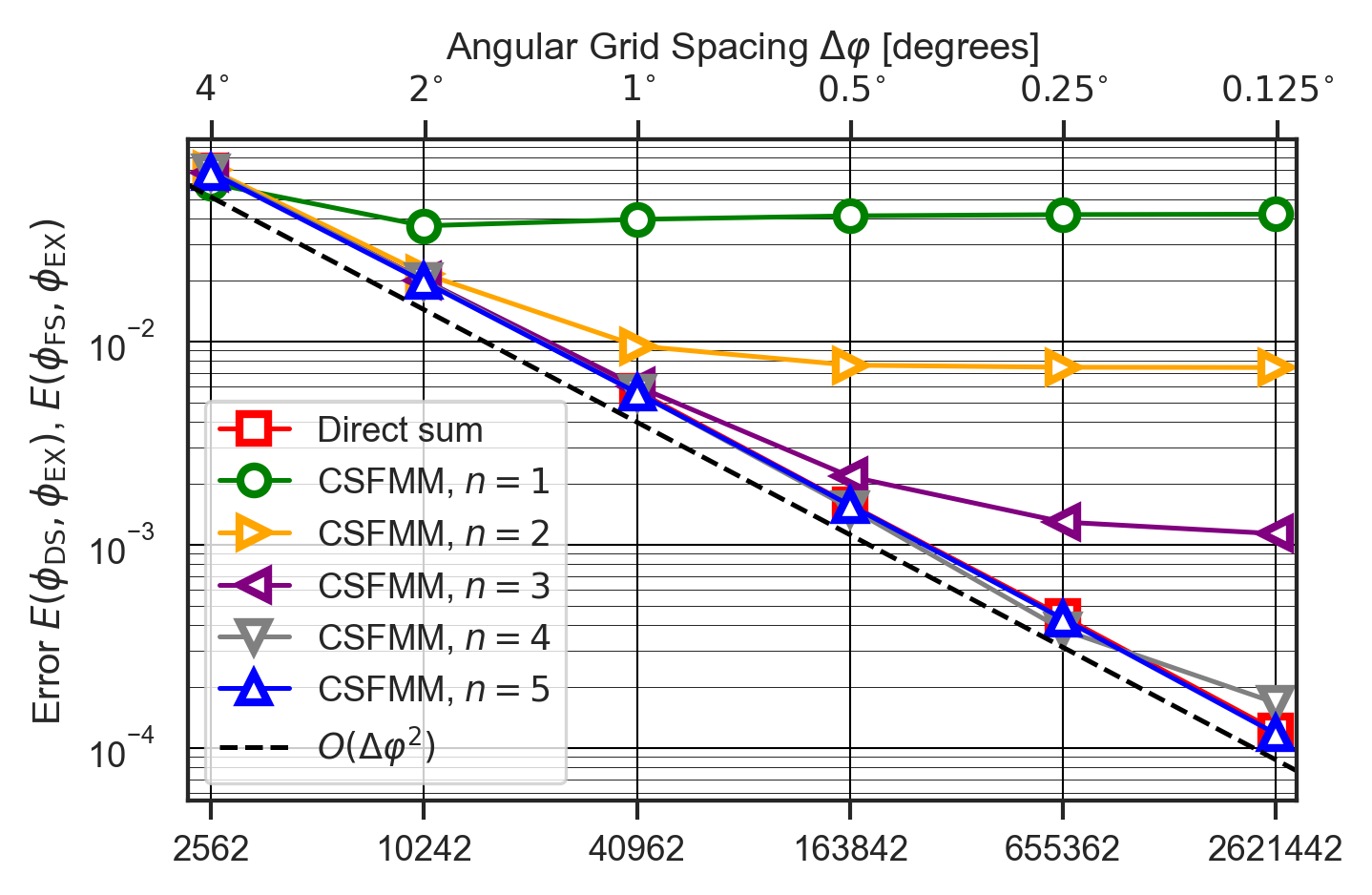}
\vskip -7.5pt
\includegraphics[width=0.80\linewidth]{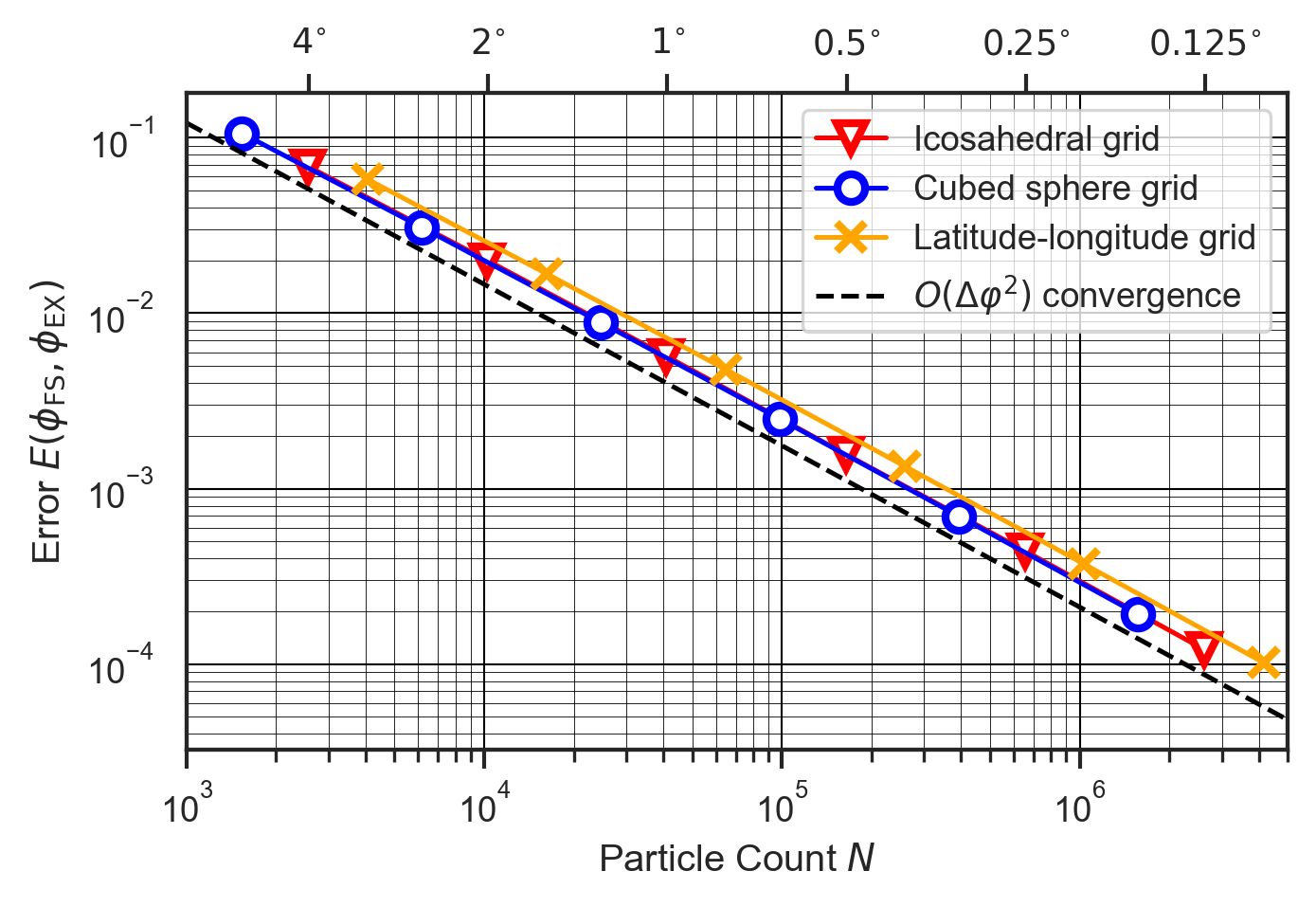}
\vskip -60pt
\setlength{\unitlength}{1cm}
\begin{picture}(5,1)
\put( -2.40, 11.50){\Large\sf{a}}
\put( -2.40,  5.20){\Large\sf{b}}
\end{picture}
\vskip 20pt
\caption{Poisson 
equation~\cref{eq:Poisson_equation},
real $Y_{4,3}$ solution,
(a) discretization error
$E(\phi_{\rm DS},\phi_{\rm EX})$
and
CSFMM solution error 
$E(\phi_{\rm FS},\phi_{\rm EX})$ 
versus particle count $N$ or 
grid spacing $\Delta\varphi$,
icosahedral grid,
CSFMM ($\mathrm{MAC}=0.7$, degree $n=1:5$),
(b) CSFMM solution error 
$E(\phi_{\rm FS},\phi_{\rm EX})$ on three grids
(icosahedral, cubed sphere,
latitude-longitude),
CSFMM ($\mathrm{MAC}=0.7$, degree $n=6$),
dashed line indicates 
2nd order convergence in $\Delta\varphi$.}
\label{fig:l2y43error}
\end{figure}

Define the global relative error,
\begin{equation}
E(\phi_1,\phi_2) =
\bigg(
\sum_{i=1}^N(\phi_1(\mathbf{x}_i)-\phi_2(\mathbf{x}_i))^2A_i\bigg/
\sum_{i=1}^N\phi_{\rm EX}(\mathbf{x}_i)^2A_i\bigg)^{\!1/2},
\end{equation} 
where
$\phi_1,\phi_2$ are chosen from 
$\phi_{\rm DS},\phi_{\rm FS},\phi_{\rm EX}$.
Consider again the 
Poisson equation~\cref{eq:Poisson_equation}
with real spherical harmonic solution $Y_{4,3}$.
\Cref{fig:l2y43error}a plots the 
discretization error 
$E(\phi_{\rm DS},\phi_{\rm EX})$
and 
CSFMM solution error
$E(\phi_{\rm FS},\phi_{\rm EX})$ 
with MAC=0.7
and
degree $n=1:5$ versus particle count $N$
or grid spacing $\Delta\varphi$
on the icosahedral grid.
The discretization error
follows the dashed line indicating
1st order convergence in $N$
(equivalent to 2nd order convergence in 
grid spacing $\Delta\varphi$).
For degree~$n=1:4$,
the CSFMM solution error follows the
dashed line for small $N$,
but it saturates at larger~$N$
when the discretization error
falls below the CSFMM approximation error;
note however that for degree~$n=5$,
the CSFMM solution error matches the 
discretization error for all~$N$. 
\Cref{fig:l2y43error}b
plots the CSFMM solution error 
$E(\phi_{\rm FS},\phi_{\rm EX})$ versus $N$
using the icosahedral, cubed sphere, 
and
latitude-longitude grids.
All three grids
exhibit 2nd order convergence in $\Delta\varphi$;
the icosahedral and 
cubed sphere grids have nearly identical error, 
while the latitude-longitude grid error is
slightly higher. 

\subsection{Error versus CSFMM parameters}

\Cref{tbl:csfmmerrors}
shows the CSFMM approximation error $E(\phi_{\mathrm{FS}},\phi_{\mathrm{DS}})$
for MAC = 0.5, 0.7, 0.9 and 
interpolation degree $n=2:2:14$.
Results are given for real spherical harmonics
$Y_{0,0}, Y_{4,3}, Y_{8,5}$,
where the integrals 
in~\cref{eq:G_L_1} and \cref{eq:G_L_2} 
were computed on the
icosahedral grid with $N=655362$
yielding discretization error
$E(\phi_{\rm{DS}},\phi_{\rm{EX}})$
between 1e-03 and 1e-04.
The CSFMM approximation error decreases 
when the MAC parameter increases or
when the interpolation degree $n$ increases,
and
in many cases it is less than the
discretization error.
For given values of MAC and degree~$n$,
the constant function $Y_{0,0}$ 
has the smallest error,
the oscillatory function $Y_{4,3}$
has larger error,
and
the more oscillatory function $Y_{8,5}$
has the largest error.

\begin{table}[htb]
\centering
\begin{tabular}{|c|c|ccccccc|}
\hline
& MAC &
$n=2$ & $n=4$ & $n=6$ & $n=8$ & $n=10$ & $n=12$ & $n=14$ \\
\hline
& 0.9 & 
6.2e-03 & 1.2e-06 & 1.2e-06 & 5.4e-08 & 3.5e-09 & 5.1e-10 & 5.9e-11 \\
$Y_{0,0}$ & 0.7 & 
3.5e-03 & 1.5e-05 & 5.5e-07 & 2.7e-08 & 5.7e-10 & 3.2e-11 & 1.7e-12 \\ 
& 0.5 & 
1.0e-03 & 5.4e-06 & 7.9e-08 & 1.3e-09 & 7.7e-12 & 2.3e-13 & 1.1e-14 \\
\hline
& 0.9 & 
1.2e-02 & 7.2e-04 & 2.2e-05 & 1.4e-06 & 6.8e-08 & 4.0e-09 & 4.5e-10 \\ 
$Y_{4,3}$ & 0.7 & 
7.5e-03 & 1.9e-04 & 5.9e-06 & 1.6e-07 & 3.0e-09 & 2.2e-10 & 8.8e-12 \\ 
& 0.5 & 
1.8e-03 & 1.2e-05 & 1.1e-07 & 1.7e-09 & 1.7e-11 & 2.3e-13 & 1.8e-14 \\
\hline
& 0.9 & 
3.3e-02 & 2.2e-03 & 9.0e-05 & 1.8e-06 & 1.2e-07 & 6.8e-09 & 7.9e-10 \\ 
$Y_{8,5}$ & 0.7 & 
1.1e-02 & 2.5e-04 & 7.7e-06 & 1.9e-07 & 5.0e-09 & 2.1e-10 & 1.6e-11 \\
& 0.5 & 
3.6e-03 & 2.5e-05 & 1.9e-07 & 2.6e-09 & 3.5e-11 & 5.3e-13 & 3.1e-14 \\
\hline
\end{tabular}
\caption{
Poisson equation~\cref{eq:Poisson_equation},
the integrals
in~\cref{eq:G_L_1} and \cref{eq:G_L_2}
are computed for 
real spherical harmonics 
$Y_{0.0}, Y_{4,3}, Y_{8,5}$,
table entries give
CSFMM approximation error $E(\phi_{\mathrm{FS}},\phi_{\mathrm{DS}})$
for MAC = 0.5, 0.7, 0.9,
interpolation degree $n=2:2:14$,
calculations used the icosahedral grid 
with $N=655362$ points yielding
discretization error
$E(\phi_{\rm{DS}},\phi_{\rm{EX}})$
between 1e-03 and 1e-04.}
\label{tbl:csfmmerrors}
\end{table}

\section{Numerical results: runtime}
\label{sec:runtime}

This section presents the runtime
for computing the $N$-body sum~\cref{eq:nbodysum}
with a given kernel and where the data field is
the real spherical harmonic $Y_{4,3}$
(the choice of data field does not 
affect the runtime).
The calculations were done on the
icosahedral grid
and
similar results hold for the
cubed sphere and latitude-longitude grids. 
All runtimes are averaged over 10 runs.

\subsection{Serial runtime versus particle count}

\begin{figure}[h!]
\centering
\includegraphics[width=0.80\linewidth]{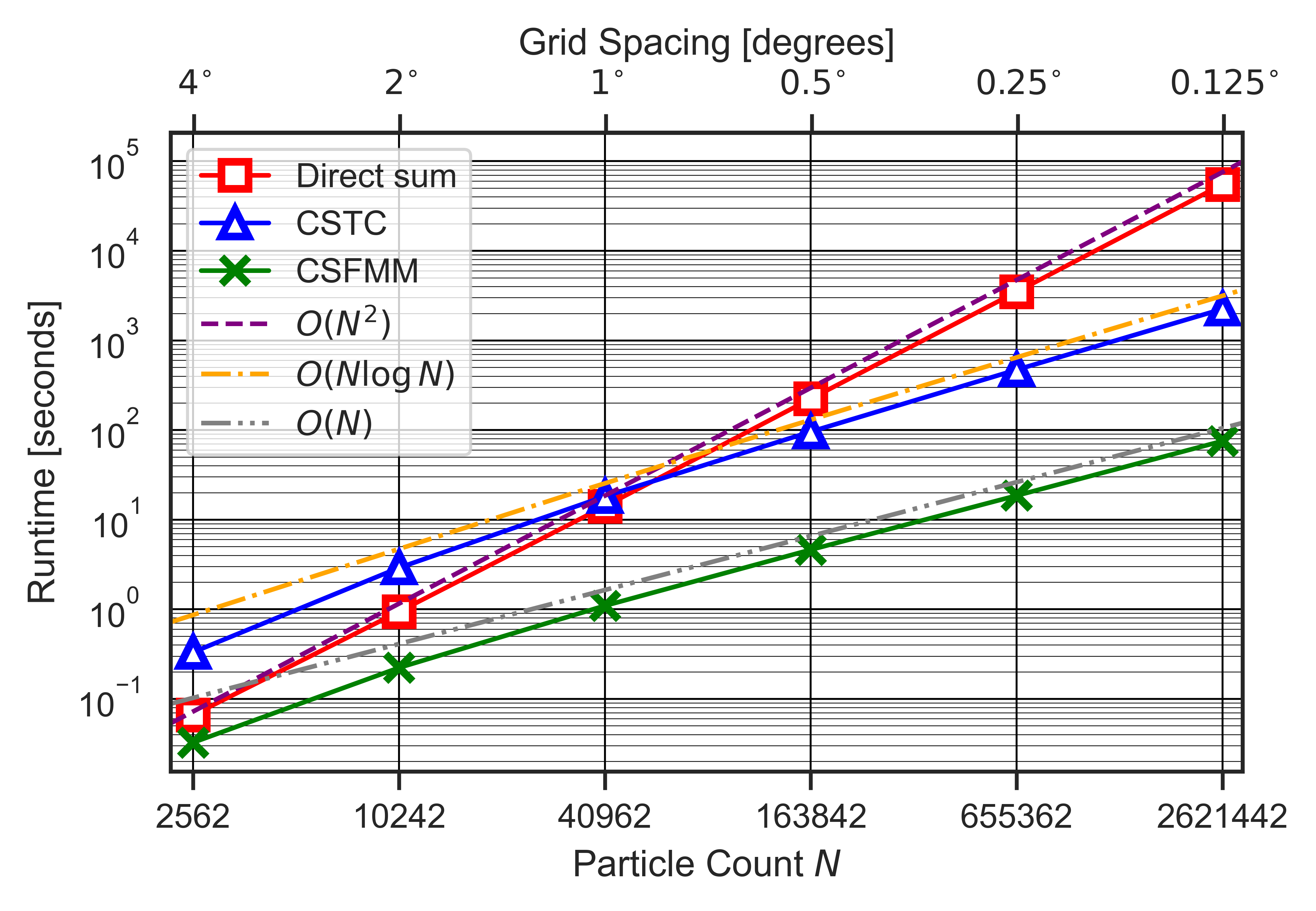}
\vskip -7.5pt
\caption{Serial runtime (s) versus 
particle count $N$ for 
$N$-body sum~\cref{eq:nbodysum}
with Laplace kernel~\cref{eq:G_L},
target and source particles on icosahedral grid,
direct sum, CSFMM and CSTC with 
$\mathrm{MAC}=0.7$, degree $n=6$.}
\label{fig:serialruntime}
\end{figure}

\Cref{fig:serialruntime} 
plots the serial runtime~(s) versus 
particle count $N$
for the Laplace kernel~\cref{eq:G_L}
with particles on the icosahedral grid.
Results are shown for direct summation, 
and
CSTC, CSFMM with MAC = 0.7, degree $n=6$
yielding fast sum approximation error 
$E(\phi_{\rm FS},\phi_{\rm DS})$ 
on the order of 1e-5.
The runtimes scale as expected,
$O(N^2)$ for direct summation, 
$O(N\log{N})$ for CSTC, 
$O(N)$ for CSFMM.  
For $N \le 40962$,
the direct sum is faster than CSTC, 
but CSTC is faster for larger $N$.
For all $N$,
CSFMM is more than an order of magnitude
faster than CSTC.
For $N=2621442$,
CSTC is 25 times faster than 
direct summation and 
CSFMM is 30 times faster than CSTC. 

\Cref{tbl:runtimes}
presents the serial runtime (s) 
versus particle count $N$ 
for the $N$-body sum~\cref{eq:nbodysum}
applied to the four kernels previously
described in \cref{sec:sphereconv}
with particles on the icosahedral grid.
The Laplace kernel~\cref{eq:G_L}
requires a logarithm, 
and 
taking it as a reference,
the biharmonic kernel~\cref{eq:G_B}
requires a dilogarithm
which is more expensive, 
the Biot-Savart kernel~\cref{eq:K_BS}
requires only simple arithmetic  
and is less expensive, 
while the SAL kernel~\cref{eq:G_SAL}
requires a logarithm
and inverse square root
and is moderately more expensive.
The runtime scales similarly with $N$
for all four kernels,
and
as $N$ increases,
the CSFMM becomes substantially faster than
direct summation.
For example with $N=655362$,
the Laplace kernel runtime is 
3490~s with direct summation
and 18.6~s with CSFMM,
while \cref{fig:l2y43error}a
shows that in this case
the CSFMM solution error 
for solving the Poisson equation
with the real $Y_{4,3}$ solution is 
about~4.4e-4. 

\begin{table}[htb]
\centering
\begin{tabular}{|c|c|c|c|c|c|c|}
\hline
$N \rightarrow$ & 2562 & 10242 & 40962 &
163842 & 655362 & 2621442 \\ 
\hline
\multicolumn{7}{|c|}{(a) direct sum runtime (s)} \\
\hline
Laplace & 0.0659 & 0.939 & 14.1 & 223 & 3490 & 55100 \\
biharmonic & 0.1030 & 1.53 & 23.1 & 418 & 7190 & 115000 \\
Biot-Savart & 0.0400 & 0.599 & 10.0 & 152 & 2700 & 54800 \\
SAL & 0.0910 &1.37 & 21.2 & 331 & 5750 & 94900 \\
\hline
\multicolumn{7}{|c|}{(b) CSFMM runtime (s)} \\
\hline
Laplace & 0.0325 & 0.222 & 1.09 & 4.61 & 18.6 & 75.4 \\
biharmonic & 0.0593 & 0.466 & 2.27 & 9.41 & 38.1 & 152 \\   
Biot-Savart & 0.0332 & 0.193 & 0.986 & 4.14 & 16.7 & 71.7 \\
SAL & 0.0514 & 0.378 & 1.85 & 7.30 & 29.4 & 121 \\
\hline
\end{tabular}
\caption{Serial runtime (s) versus particle count $N$
for $N$-body sum~\cref{eq:nbodysum},
icosahedral grid,
four kernels,
Laplace~\cref{eq:G_L}, biharmonic~\cref{eq:G_B}, 
Biot-Savart~\cref{eq:K_BS}, 
SAL~\cref{eq:G_SAL},
(a) direct summation,
(b) CSFMM ($\mathrm{MAC}=0.7$, degree $n=6$).}
\label{tbl:runtimes}
\end{table}

\subsection{Serial runtime versus error}
 
\begin{figure}[h!]
\centering
\includegraphics[width=\textwidth]{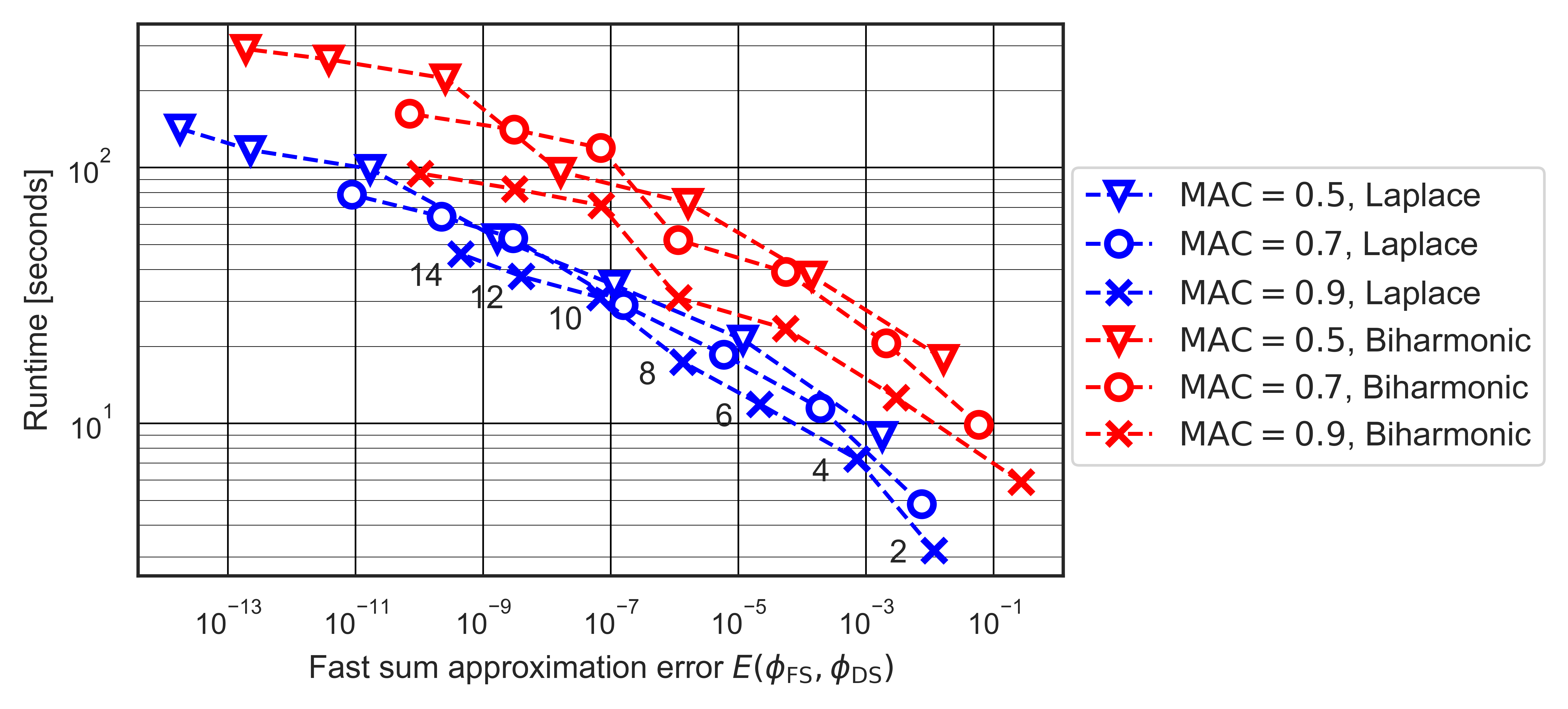}
\vskip -10pt
\caption{Serial runtime (s)
versus CSFMM approximation error 
$E(\phi_{\rm FS},\phi_{\rm DS})$
for computing $N$-body sum~\cref{eq:nbodysum},
Laplace kernel~\cref{eq:G_L},
biharmonic kernel~\cref{eq:G_B},
icosahedral grid, particle count $N=655362$,
CSFMM with MAC = 0.5, 0.7, 0.9, 
degree $n=2:2:14$ 
increasing from lower right to upper left
on each connected line.}
\label{fig:varydtheta}
\end{figure}

\Cref{fig:varydtheta} 
plots the serial runtime (s) versus
CSFMM approximation error
$E(\phi_{\rm FS},\phi_{\rm DS})$
for computing the $N$-body sum~\cref{eq:nbodysum}
with the
Laplace kernel~\cref{eq:G_L}
and
biharmonic kernel~\cref{eq:G_B}.
The calculation used the
icosahedral grid with particle count $N=655\,362$.
Results are shown for CSFMM
with MAC = 0.5, 0.7, 0.9
and degree $n=2:2:14$ increasing
from lower right to upper left
on each connected line.
The error decreases 
for smaller MAC and larger degree $n$,
reaching close to machine precision
for MAC = 0.5 and degree $n=14$. 
The runtime increases as the error decreases.
The biharmonic runtime is slightly more
than twice the Laplace runtime,
which is consistent with the 
results in~\Cref{tbl:runtimes}
using MAC = 0.7 and $n=6$.
The results provide guidance on the best choice
of MAC parameter and interpolation degree~$n$
to achieve a given error with the
least runtime.

\subsection{Parallel runtime}

\begin{figure}[htb]
\centering
\includegraphics[width=0.9\textwidth]{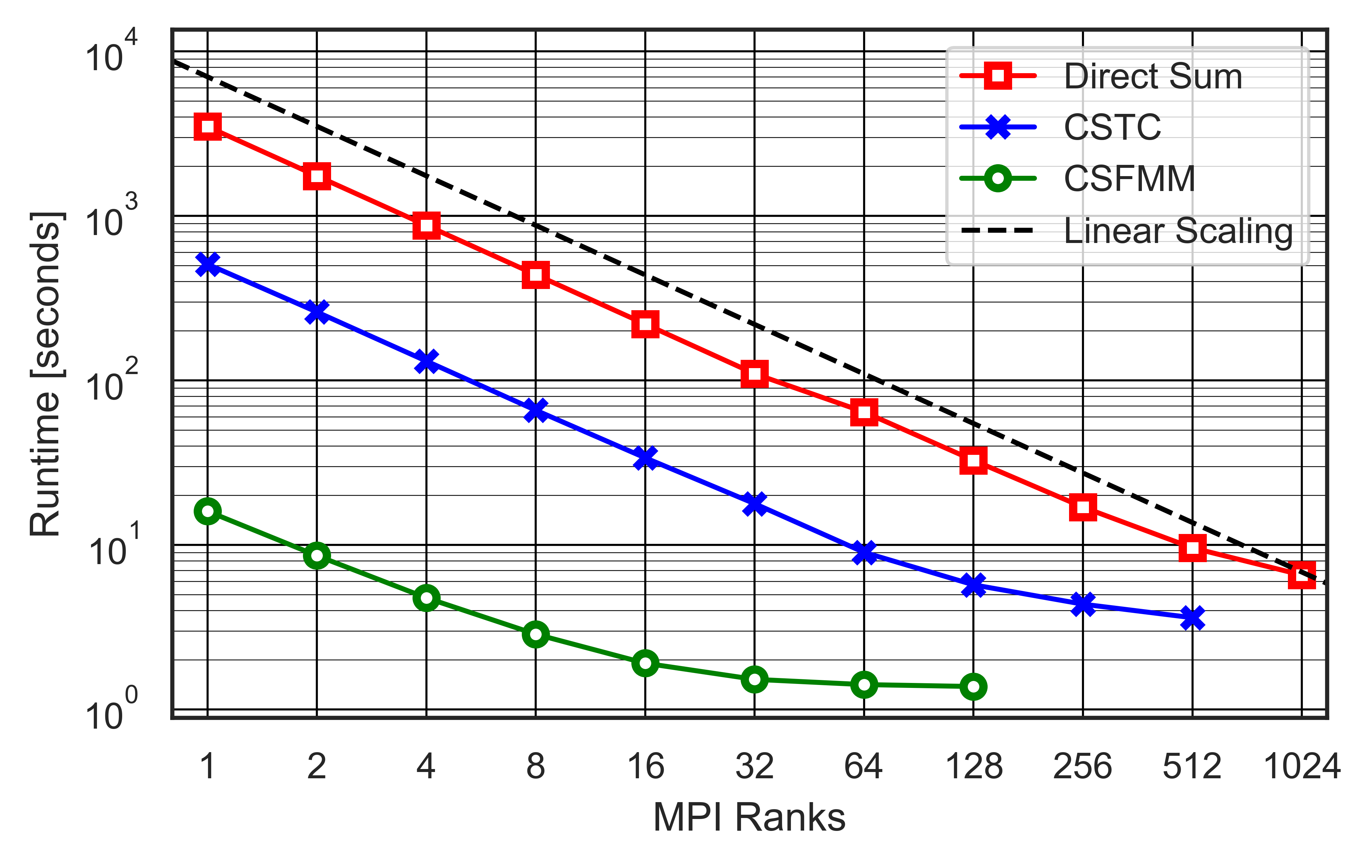}
\vskip -10pt
\caption{Parallel runtime (s) versus 
number of MPI ranks for
computing $N$-body sum~\cref{eq:nbodysum}
with Laplace kernel~\cref{eq:G_L},
icosahedral grid, $N=655362$,
results for
direct sum, 
CSTC and CSFMM with 
$\mathrm{MAC}=0.7$, degree $n=6$.}
\label{fig:parallelruntime}
\end{figure}

The direct sum, CSFMM, and CSTC were parallelized 
using MPI for 
intra-node and inter-node parallelism,
with communication by 
Remote Memory Access~\cite{mpirma}. 
In this work each MPI rank
computes an equal portion of the
PP, PC, CP, and CC interactions,
followed by a global reduction. 
The upward and downward passes in the
CSFMM and CSTC were computed serially  
since they consume only
a small fraction of the runtime.
\Cref{fig:parallelruntime}
plots the parallel runtime (s)
versus number of MPI ranks
for computing the $N$-body sum~\cref{eq:nbodysum}
with the Laplace kernel~\cref{eq:G_L}
on the icosahedral grid
with particle count $N=655362$.
Results are shown for direct summation,
and CSTC and CSFMM with
$\mathrm{MAC}=0.7$ and degree $n = 6$.
All three methods display near-linear scaling 
for small processor counts. 
The parallel efficiency falls to around 50\%
at $1024$ ranks for direct summation 
(runtime 7~s), 
$256$ ranks for CSTC (runtime 4.5~s),
and
32~ranks for CSFMM (runtime 1.5~s). 
The CSFMM runtime saturates
primarily due to the serial computation of the 
upward and downward passes.
The next two sections demonstrate
the capability of the CSFMM for problems
in geophysical fluid dynamics.

\section{Barotropic vorticity equation}
\label{sec:BVEres}

This section considers several examples 
in which the
barotropic vorticity equation (BVE)
on a rotating sphere~\cref{eq:BVE} 
is solved by a remeshed vortex 
method~\cite{bosler2014lagrangian,koumoutsakos1997inviscid,chorin1973numerical}.
The sphere rotation rate is
$\Omega = 2\pi\,\mathrm{ day}^{-1}$.
The particles ${\bf x}_i(t)$ are
initialized to lie at
icosahedral grid cell centers~(\cref{fig:grids}a).
The ODEs~\cref{eq:BVE_ODEs} are solved 
by the 4th order Runge-Kutta method with
time step $\Delta t=0.01$~day
and
remeshing is done at every time step
by local biquadratic interpolation.
The velocity~\cref{eq:K_BS}
is computed by direct summation or
CSFMM (MAC = 0.7, degree $n=6$) as indicated.
The chosen parameters ensure that the
time stepping error 
and
CSFMM approximation error
are smaller than
the midpoint rule discretization error.

\begin{figure}[htb]
\centering 
\includegraphics[width=\textwidth]{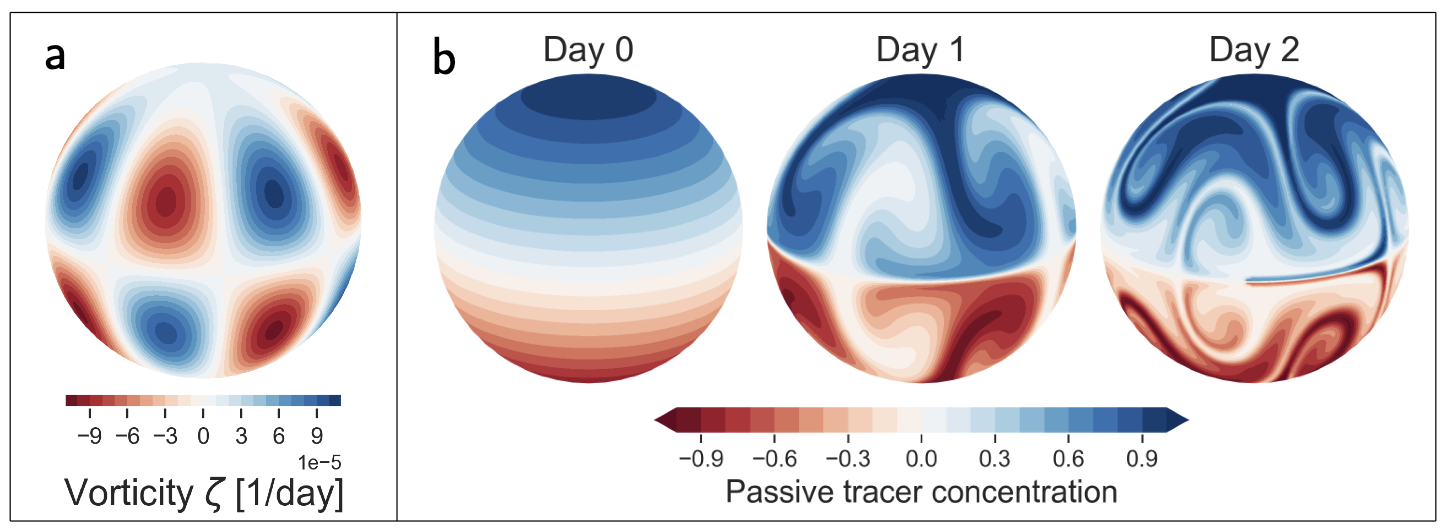}
\vskip -5pt
\caption{
Barotropic vorticity equation~\cref{eq:BVE},
(a) RH wave vorticity~\cref{eq:RH_initial},
(b) passive tracer advection by RH wave,
tracer field at Day 0, 1, 2, 
icosahedral grid, 
$N=163842$ particles,
velocity~\cref{eq:K_BS} 
computed by CSFMM (MAC=0.7, $n=6$),
a video is available in the 
Supplementary Material.}
\label{fig:bvepassivetra}
\end{figure}

\begin{figure}[htb]
\centering
\includegraphics[width=0.8\linewidth]{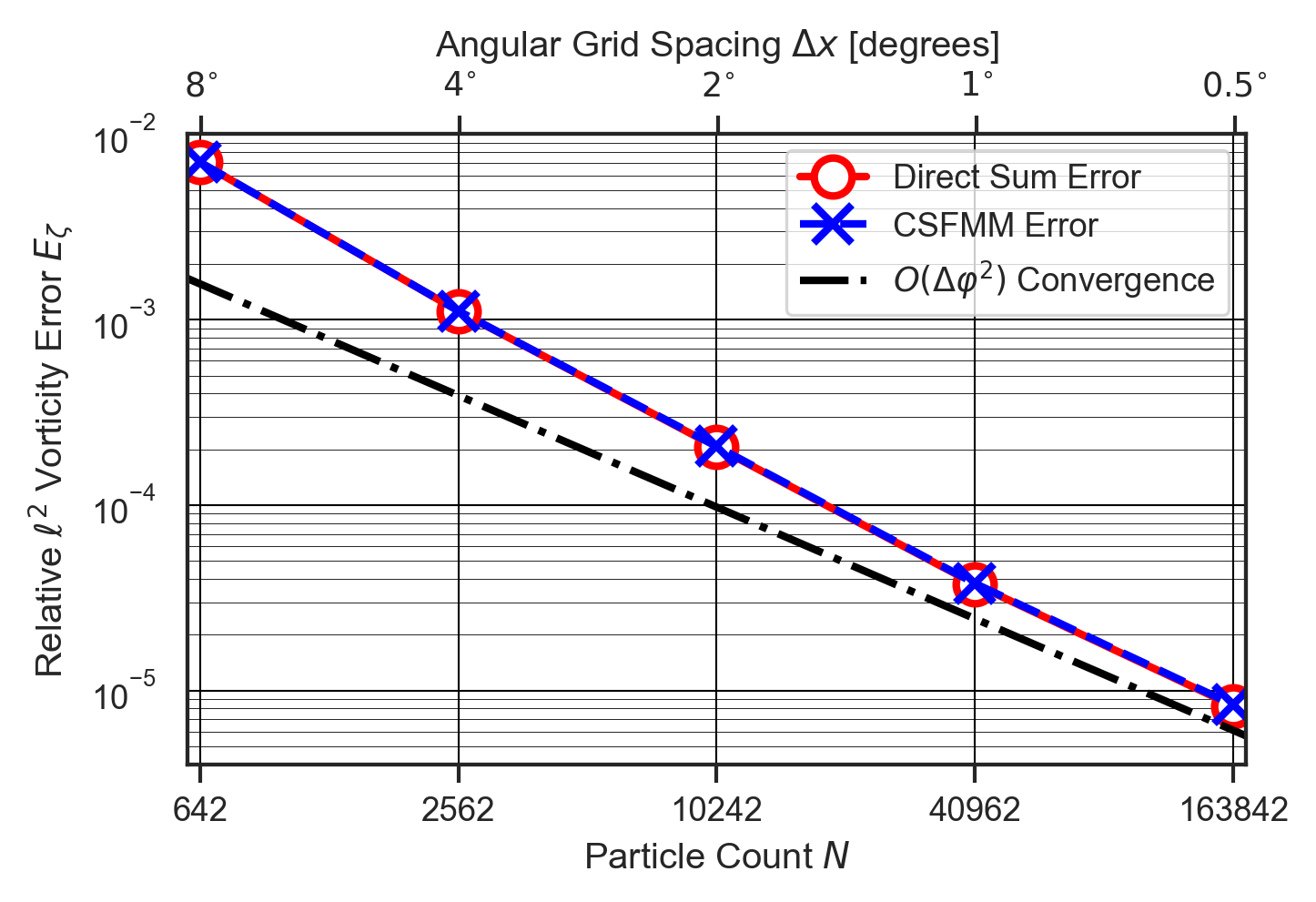}
\vskip -5pt
\caption{
Barotropic vorticity equation~\cref{eq:BVE},
RH wave vorticity~\cref{eq:RH_initial},
vorticity error 
$E_{\zeta}$ at time $t=1$~day~\cref{eq:bveerr}
versus particle count~$N$,
icosahedral grid,
particle velocity~\cref{eq:K_BS}
computed by direct sum and 
CSFMM (MAC=0.7, degree $n=6$).}
\label{fig:bveerr}
\end{figure}

\vskip 5pt
\noindent
{\bf Rossby-Haurwitz wave.}
Rossby-Haurwitz (RH) waves are a class of 
traveling wave solutions of the BVE
on a rotating sphere~\cite{haurwitz1940motion}.
We choose a wave with 
initial vorticity
\begin{equation}
\label{eq:RH_initial}
\zeta_0(\theta,\lambda) =
\frac{2\pi}{7}\sin\theta +
30\sin\theta\cos^4\!\theta\cos 4\lambda, 
\end{equation} 
where the first term is
a solid body rotation of the sphere 
and
the second term is a real spherical harmonic
with degree $n=5$ 
and wavenumber $m=4$.
\Cref{fig:bvepassivetra}a
shows the vorticity $\zeta_0$ 
which is anti-symmetric about the equator
with alternating
positive and negative vortex cores
in each hemisphere.
\Cref{fig:bvepassivetra}b
shows a tracer advection calculation explained below.
The solid body rotation
ensures that the RH wave speed vanishes, 
so the vorticity is stationary, 
$\zeta({\bf x},t) = \zeta_0({\bf x})$,
and this provides a reference
for error assessment;
note however that the fluid velocity is nonzero, 
so the computational particles~${\bf x}_i(t)$
are moving and
their vorticity~$\zeta_i(t)$ is changing.
The error in the computed vorticity at time $t$ is defined by
\begin{equation}\label{eq:bveerr}
E_{\zeta}(t) =
\left(
\sum_{i=1}^N
(\zeta_i(t) - \zeta_0(\mathbf{x}_i))^2A_i
\bigg/\sum_{i=1}^N\zeta_0(\mathbf{x}_i)^2A_i
\right)^{\!1/2}.
\end{equation}
\Cref{fig:bveerr} plots $E_\zeta$ 
at time $t=1$ day
versus particle count $N$,
where the particle velocity~\cref{eq:K_BS}
was computed two ways,
by direct summation and CSFMM.
Both methods yield essentially 
the same error,
showing that the CSFMM approximation error
is negligible compared to the
discretization error.
As $N$ increases,
the error asymptotes to the 
dashed-dot line indicating 
convergence at the rate $O(\Delta\varphi^2)$.
For the finest grid with $N=163\,842$,
the CSFMM is 36 times faster than
direct summation.

\vskip 5pt
\noindent
{\bf Passive tracer advection by a RH wave.}
The transport of a passive tracer is a
common test case for the performance of
numerical advection schemes on the sphere,
where the tracer field $q({\bf x},t)$
represents the concentration of a 
nutrient, pollutant,
or other substance~\cite{bosler2017lagrangian,
bradley2024stabilized,
kent2014dynamical,
lauritzen2014standard}.
Here we consider passive tracer
advection in the velocity field
${\bf u}({\bf x},t)$
defined by the RH wave
with vorticity~\cref{eq:RH_initial}
as depicted in \Cref{fig:bvepassivetra}a.
The governing equation is 
$\partial q/\partial t + 
{\bf u} \cdot \nabla q = 0$,
where the velocity is computed by
the vortex method and CSFMM.
In this case the 
computational particles ${\bf x}_i(t)$
carry tracer values $q_i(t)$
in addition to vorticity~$\zeta_i(t)$~\cite{bosler2014lagrangian}. 

\Cref{fig:bvepassivetra}b
shows results where 
the initial tracer field 
is the particle $z$-coordinate;
a video is available in the 
Supplementary Material.
On Day 0,
the tracer varies smoothly in latitude
and is independent of longitude,
with a band of low values near the equator
and 
high value caps at the poles.
Later on the tracer field becomes more complex, 
but still anti-symmetric about the equator,
so we will describe the dynamics
in the northern hemisphere.
On Day~1,
the low value band initially
near the equator has rolled up around the
clockwise-rotating negative vortex cores,
while the 
counterclockwise-rotating positive vortex cores
entrain high value tracer from the north pole.
On Day~2,
the tracer field has a staggered pattern,
with low value clockwise spiral filaments
at a low latitude
and
high value counterclockwise spiral filaments
at a higher latitude. 

Note that the tracer field in this problem
has certain properties;
for example,
the total tracer mass vanishes by symmetry,
$\int_S q({\bf x},t)dS = 0$,
and
the tracer field is bounded,
$|q({\bf x},t)| \le 1$.
If desired,
these properties can be enforced 
to machine precision 
in the numerical solution using several techniques
(this is called {\it property preservation},
e.g.~see~\cite{bosler2019conservative,
bradley2021islet,
bradley2019communication}
and the references therein),
but this was not implemented
in the present calculations;
instead the numerical parameters were chosen
to ensure that the property errors have
negligible effect.
To illustrate,
consider computed values of the
tracer mass~$m$
and tracer bound~$q_{\rm max}$ defined by
\begin{equation}
m = \sum_{i=1}^N q_iA_i,
\quad
q_{\rm max} = \max\{|q_i|, i=1:N\},
\end{equation}
where $q_i$ 
is the computed tracer value at 
particle ${\bf x}_i$ 
at the end of the calculation on Day~2.
\Cref{tab:tracer_errors}
shows the effect of the point count $N$,
where the errors in the
tracer mass and tracer bound 
diminish as the spatial grid is refined.
With $N = 163842$
the errors are too small to be visible
in \cref{fig:bvepassivetra}b.

\begin{table}[htb]
\centering
\begin{tabular}{|c|c|c|c|}
\hline
tracer error &
$N = 10242$ & $N = 40962$ & $N = 163842$ \\
\hline
$|m|$ & 
6.49e-7 & 2.56e-7 & 1.83e-7 \\
$q_{\rm max}-1$ &
2.36e-2 & 5.26e-3 & 1.29e-3 \\
\hline
\end{tabular}
\caption{
Tracer property errors versus point count $N$,
mass error $|m|$, 
bound error $q_{\rm max}-1$,
errors are reported on Day~2 in 
\cref{fig:bvepassivetra}b,
tracer advected in RH wave with
velocity computed by vortex method
on isocahedral grid,
CSFMM (MAC=0.7, degree $n=6$).}
\label{tab:tracer_errors}
\end{table}

\begin{figure}[htb]
\centering
\includegraphics[width=\linewidth]{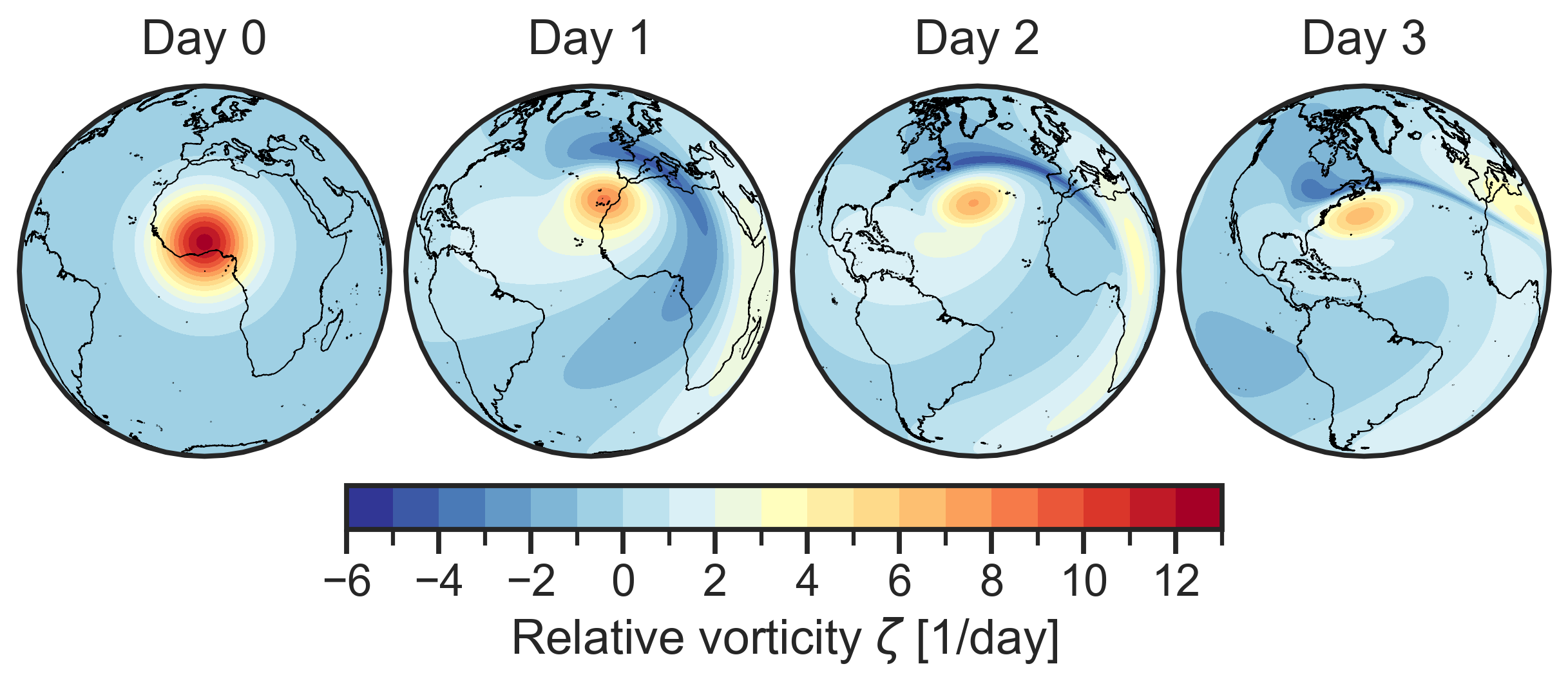}
\vskip -10pt
\caption{
Barotropic vorticity equation~\cref{eq:BVE}, 
Gaussian vortex,
Day 0:~initial vorticity~\cref{eq:Gaussian_ic}, 
Day 1,\,2,\,3:~vorticity computed by a
remeshed vortex method with
CSFMM ($\mathrm{MAC}=0.7$, degree $n=6$),
iscosahedral grid,
$N = 163\,842$ particles,
videos in the 
Supplementary Material show the evolution
of the vorticity 
and a passive tracer.}
\label{fig:gv}
\end{figure}

\vskip 5pt
\noindent
{\bf Gaussian vortex.}
\Cref{fig:gv} shows the evolution of a
Gaussian vortex as a crude model of a cyclone
in the northern hemisphere.
The initial vorticity,
\begin{equation}
\label{eq:Gaussian_ic}
\zeta_0(\theta,\lambda) =
4\pi\exp(-16|\mathbf{x}-\mathbf{x}_c|^2) -
0.196353,
\quad \mathbf{x}_c = (\pi/20,0),
\end{equation}
has a positive vortex core
slightly above the equator
and a negative constant background
to ensure zero total 
vorticity~\cite{bosler2014lagrangian}.
The continents are shown to help
visualize the vortex motion;
the BVE has no topography.
The calculation used $N=163\,842$ particles
on the icosahedral grid with
velocity~\cref{eq:K_BS} computed by CSFMM.

Recall that the BVE~\cref{eq:BVE}
conserves the absolute vorticity $\zeta+f$ 
of fluid particles;
hence if a particle in the northern hemisphere
moves to higher latitude,
its planetary vorticity $f$ increases
and its relative vorticity $\zeta$ decreases;
conversely, if such a particle moves to lower latitude, $\zeta$ increases.
On Day~0,
the fluid particles start rotating
counterclockwise around the vortex core,
and as their latitude changes,
the vorticity above the core decreases
and
below the core it increases;
hence the vorticity field acquires a
dipole component that advects the
core to the northwest.
On Day~1,
the vortex core has weakened slightly
and
a layer of negative vorticity has formed
on its northeast side 
stretching around the core
and reaching back close to the 
initial core location.
The negative vortex layer advects
fluid particles behind it to lower latitude and
this creates a weak positive vortex layer
further behind.
On Day~2,
the vortex core is propagating westward
with a streamwise-aligned elliptical shape,
while the negative vortex layer 
has become thinner.
On Day~3,
the vortex core has weakened
to about half its initial amplitude,
and
the head of the negative vortex layer
is flattening into a core;
however, the positive core is stronger 
than the negative core,
and the dipole path bends to the southwest.
Note that as the sphere completes
three rotations,
the vortex travels from West Africa to the 
east coast of North America, 
leaving behind several 
alternating vortex layers in its wake.
Videos in the 
Supplementary Material show the
evolution of the vorticity 
and
a passive tracer.

\section{Self-attraction and loading}
\label{sec:SAL}

\begin{figure}[htb]
\includegraphics[width=\linewidth]{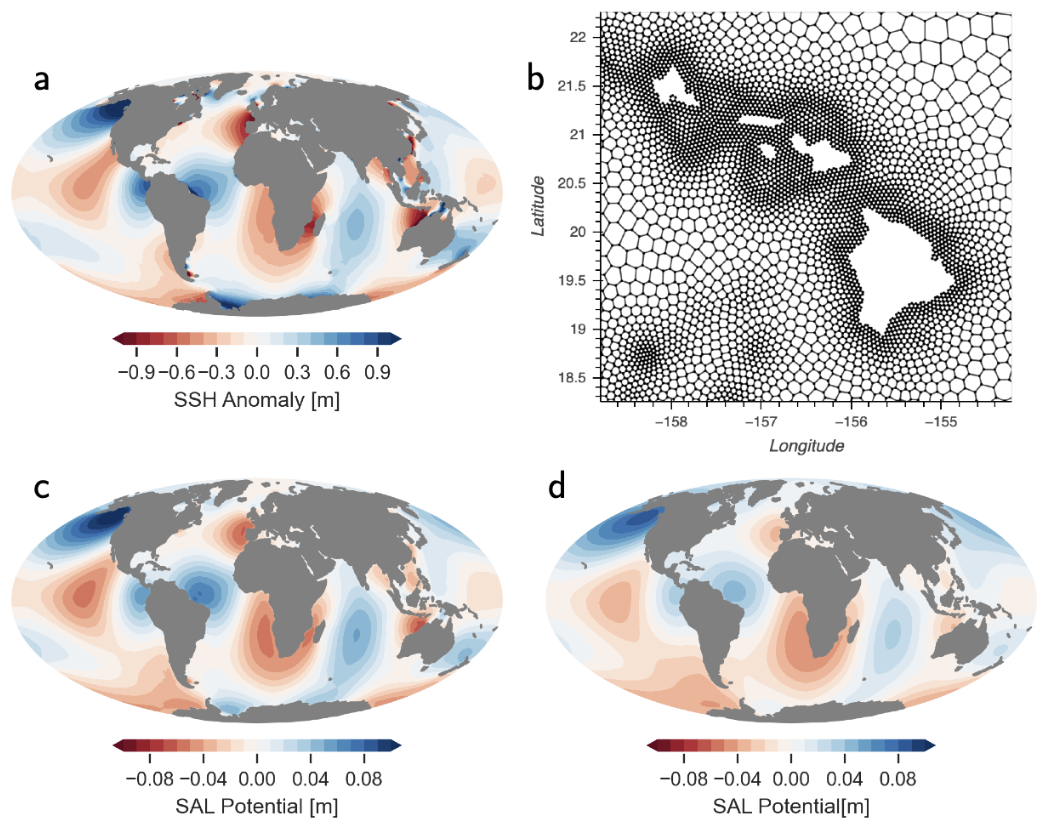}
\vskip -5pt
\caption{SAL calculation,
(a) input sea surface height anomaly~$\eta$ 
from MPAS-Ocean calculation with 
$N=2313486$ cells,
(b) detail of variable resolution hexagonal
grid near Hawaii,
SAL potential~$\eta_{\rm SAL}$ computed by
(c) spherical harmonic series, 
(d) convolution/CSFMM 
(MAC = 0.7, degree $n=2$).}
\label{fig:sal}
\end{figure}

The next example is a calculation
of the SAL potential $\eta_{\rm SAL}$ 
as discussed in \Cref{sec:sal}.
The input field is the
sea surface height anomaly~$\eta$ 
shown in~\cref{fig:sal}a
which was computed using MPAS-Ocean
(Model for Prediction Across Scales~\cite{ringler2013multi})
on a variable resolution hexagonal grid with 
$N=2313486$ cells.
The grid spacing varies 
from 45 km in the open ocean to 
5 km in coastal regions
and
\cref{fig:sal}b shows a detail of the
grid near Hawaii.

\Cref{fig:sal}
also shows the SAL potential $\eta_{\rm SAL}$
computed two ways,
(c) spherical harmonics series~\cref{eq:SH}
using MPAS-Ocean~\cite{barton2022global,
brus2023scalable},
(d) spherical convolution~\cref{eq:spherical_convolution} 
of the SAL 
kernel $G_{\rm SAL}$~\cref{eq:G_SAL_approx}
with the sea surface height anomaly~$\eta$
on the MPAS-Ocean grid using 
CSFMM (MAC=0.7, degree $n=2$) to accelerate the 
$N$-body sum~\cref{eq:nbodysum}. 
The two computed SAL potentials are similar
and both resemble the
sea surface height anomaly $\eta$ 
in~\cref{fig:sal}a
scaled by a factor of approximately 1/10
(this scaling has been noted previously~\cite{accad1978solution}).
However, close examination 
reveals some differences;
the convolution result has
lower amplitude and is smoother than the
spherical harmonics result
especially in some coastal regions.
This is important because
tidal calculations require
the SAL acceleration $-g\nabla\eta_{\rm SAL}$,
and
small errors in the potential~$\eta_{\rm SAL}$
will be amplified by differentiation.
Hence, while the results 
support the validity of the
convolution approach,
higher resolution calculations
are needed to determine which method is more accurate.
In terms of efficiency,
the convolution/CSFMM runtime was 
approximately 90 seconds in serial,
however the spherical harmonics result
came from a multicore MPAS-Ocean simulation 
that computed more than $\eta_{\rm SAL}$
and the
relevant runtime is not available.
Further comparison of the two approaches
is planned once the 
convolution/CSFMM scheme is installed in 
an ocean general circulation model.

\section{Summary}
\label{sec:summary}

This work introduced the
Cubed Sphere Fast Multiple Method (CSFMM)
for summing pairwise particle 
interactions~\cref{eq:nbodysum}
that arise from discretization of
integral transforms and
convolutions on the sphere.
The kernel approximations use
barycentric Lagrange interpolation on
a quadtree composed of cubed sphere grid cells.
The scheme is kernel-independent
and
requires kernel evaluations 
only at points on the sphere.
The midpoint rule was used to
discretize the integral transforms
with three spherical grids
(icosahedral, cubed sphere, 
latitude-longitude)
and
2nd order convergence in the
angular grid spacing $\Delta\varphi$
was observed
(equivalent to
1st order convergence in the 
point count $N$).
A Cubed Sphere Tree Code (CSTC)
was described for comparison,
and
for a given level of accuracy,
the CSTC runtime scales like $O(N\log N)$,
while the CSFMM runtime scales like $O(N)$.
In parallel calculations,
both schemes exhibit near-linear speedup
until the parallel runtime falls to a 
small fraction of the serial runtime.
Results were presented for the
Poisson and biharmonic equations on the sphere,
barotropic vorticity equation on a rotating sphere,
and the
self-attraction and loading potential in
tidal calculations,
where in the latter case
the CSFMM was applied to particles on a 
variable resolution hexagonal grid.
Note that the cubed sphere is used here
in two different ways, 
(1) it is one of the three 
spherical partitions shown
in~\Cref{fig:grids}
used to discretize the convolution integral,
(2) the cubed sphere grid cells provide
the particle clusters used in the 
CSFMM and CSTC to 
accelerate the calculation of the $N$-body sum.

There are several directions for future study. 
In terms of code development, 
the upward and downward passes 
can be parallelized, 
and we plan to develop a 
mixed MPI/Kokkos implementation of CSFMM
for improved portability and multinode
GPU capability~\cite{edwards2014kokkos}.
In terms of applications,
we plan to install the 
convolution/CSFMM SAL scheme as an option
in MOM6 (Modular Ocean 
Model~6~\cite{wang2024improving}),
which currently uses a spherical harmonics 
SAL scheme.
The SAL acceleration is also needed
in regional tide 
prediction~\cite{irazoqui2017effects}, 
where the global nature of 
spherical harmonics limits their capability,
and we plan to explore
the convolution/CSFMM approach in this case as well.
Finally, another goal is extending the 
vortex method 
to the shallow water equations on the sphere,
with particles that carry vorticity and
divergence,
using discrete convolutions accelerated by 
the CSFMM to compute the velocity,
and tree-based
adaptive mesh refinement of the vorticity and divergence fields~\cite{bosler2014lagrangian,bosler2017lagrangian,
sandberg2025farsight,xu2023dynamics}.

\section*{Acknowledgments}
The authors thank Brian Arbic,
Peter Bosler and Christiane Jablonowski   
for helpful discussions,
and Kristin Barton 
for providing the MPAS data. 
We gratefully acknowledge high-performance computing resources on the Derecho system provided by the NSF-sponsored
National Center for Atmospheric Research (NCAR).
This work was also supported in part through resources provided by Advanced Research Computing at the University of Michigan. 

\bibliographystyle{siamplain}
\bibliography{references}
\end{document}